\documentclass{article}
\usepackage{amsfonts}
\newcommand{\Ni}{\mathbb{N}}

\newcommand{\KK}{\mathbb{K}}
\newcommand{\CC}{\mathbb{C}}
\newcommand{\NN}{\mathbb{N}}
\newcommand{\RR}{\mathbb{R}}
\newcommand{\QQ}{\mathbb{Q}}

\newcommand{\UU}{\mathbb{U}}
\newcommand{\ZZ}{\mathbb{Z}}\newcommand{\notequiv}{\equiv \kern-0.95em / \mbox{ }}
\newcommand{\QED}{{\unskip\nobreak\hfil\penalty50
\hskip1em\hbox{}\nobreak\hfil $\Box$
\parfillskip=0pt \finalhyphendemerits=0 \par\medskip\noindent}}
\newtheorem{num}{}[section]
\newtheorem{defin}[num]{}
\newenvironment{dem}{}{\QED}
\newenvironment{defi}{\begin{defin} \rm}{\end{defin}}
\newenvironment{address}{\begin{center}\small\rm}{\end{center}}
\title{c-regular cyclically ordered groups.}
\author{G. Leloup \& F. Lucas}
\date{December 17, 2013}
\begin{document}
\maketitle
\begin{abstract} 
We define and we characterize regular and c-regular cyclically ordered abelian groups. 
We prove that every dense c-regular cyclically ordered abelian group is elementarily equivalent to 
some cyclically ordered group of 
unimodular complex numbers, that every discrete c-regular cyclically ordered abelian 
group is elementarily equivalent to 
some ultraproduct of finite cyclic groups, and that the discrete regular non-c-regular 
cyclically ordered abelian groups 
are elementarily equivalent to the linearly cyclically ordered group $\ZZ$. 
\end{abstract}
\section{Introduction.}
Unless otherwise stated the groups are abelian groups. 
The definitions and basic properties of cyclically ordered groups can be found 
for example in \cite{Fu 63}, \cite{GLL 12} or \cite{Le 12}. For the reader's convenience 
we recall them in Section \ref{par2}. 
Being c-archimedean (i.e. isomorphic to a subgroup of the cyclically ordered group 
$\KK$ of the complex 
numbers of module $1$) for a cyclically ordered group is not a first order property. 
In this paper we define the c-regular cyclically ordered groups, and we 
prove that the class of all c-regular cyclically ordered groups is 
the lowest elementary class containing all 
the c-archimedean cyclically ordered groups. 
We also define the regular cyclically ordered 
groups, and we prove that the class of all regular cyclically ordered groups is 
the union of the class of all c-regular 
cyclically ordered groups and of the class of all linearly cyclically ordered groups.  \\ 
\indent We prove that an abelian cyclically ordered group $G$ which is not 
c-ar\-chi\-me\-dean is 
c-regular if and only if its linear part $l(G)$ is a regular linearly ordered group and the 
cyclically ordered group $K(G)=(G/l(G))$ is divisible and its torsion group is 
isomorphic to the group of all roots of $1$ in the field of all complex numbers. Another cha\-rac\-te\-ri\-za\-tion is that its unwound is a regular linearly ordered group. \\ 
\indent Every regular cyclically ordered group is either dense or discrete. We prove that any two 
dense c-regular cyclically ordered groups are elementarily equivalent if and only if they are 
elementarily equivalent in the language of groups, and that any dense c-regular group is 
elementarily equivalent to some c-archimedean group. 
In the case of discrete groups, we prove that 
every discrete c-regular group is e\-le\-men\-ta\-ri\-ly 
equivalent to an ultraproduct of finite cyclic groups. 
Every discrete c-regular group $G$ which is not c-archimedean contains an 
elementary substructure $C$ such that $K(C)\simeq \UU$ and $l(C)\simeq \ZZ$. 
We define first order formulas $D_{p^n,k}$ ($p$ a prime number, $n\in \NN$ and 
$k\in \{0,\dots, p^n-1\}$) such that any two discrete non-c-archimedean c-regular groups 
$G$ and $G'$ are elementarily equivalent if and only if $G$ and $G'$ satisfy the same formulas $D_{p^n,k}$. 
If $C$ and $C'$ are the elementary substructures of $G$ and $G'$, respectively, 
which we defined above, then $G\equiv G'$ if and only if $C\simeq C'$. By means of 
the formulas $D_{p^n,k}$ which are satisfied by $G$, we define a family of subsets 
of $\NN^*$, which we will call the family of subsets of $\NN^*$ characteristic of 
$G$, such that for every non principal ultrafilter $U$ on $\NN^*$, $G$ is 
elementarily equivalent to the ultraproduct of the cyclically ordered groups  
$\ZZ/n\ZZ$ modulo $U$ if and only if $U$ contains the family of subsets 
of $\NN^*$ characteristic of $G$. Furthermore, for every non-c-archimedean 
c-regular cyclically ordered group, there exists an ultrafilter containing the family of subsets 
of $\NN^*$ characteristic of this group. 
Remark that the class of all discrete regular cyclically ordered groups is the lowest 
elementary class 
which contains all the cyclically ordered groups generated by one element. Furthermore, 
this class is divided into two disjoint elementary subclasses, the class of discrete 
c-regulars cyclically ordered groups, and the class of discrete regular linearly cyclically 
ordered groups, that is, the groups which are elementarily equivalent to the linearly cyclically ordered group $\ZZ$. \\[2mm]  
\indent In the last section, we deal with the same groups, wihtout a predicate for the 
cyclic order, to show which informations this predicate carries. \\
\section{Cyclically ordered groups.}\label{par2} 
\indent Recall that a {\it cyclically 
ordered group} is a group $G$ (which is not necessarily abelian) together with a ternary 
relation $R$ which satisfies:  \\
(1) $\forall (a,b,c)\in G^3$, 
$R(a,b,c) \Rightarrow a \neq b \neq c \neq a$ ($R$ is ``strict''), \\
(2) $\forall (a,b,c)\in G^3$, $R(a,b,c) \Rightarrow R(b,c,a)$ ($R$ is cyclic), \\
(3) for all $c \in G$, $R(c , \cdot , \cdot )$ defines a linear order on  $G \backslash \{ c \}$, \\ 
(4) $\forall (a,b,c,u,v) \in G^5, \; R(a,b,c) \Rightarrow R(uav,ubv,ucv)$ 
($R$ is compatible). \\
\indent 
For example, let $\KK$ be the group of all complex numbers of module $1$, for $e^{i\theta_1}$, 
$e^{i\theta_2}$, $e^{i\theta_3}$ in $\KK$ with $0\leq \theta_i<2\pi$, ($i\in \{1,2,3\}$) 
we let $R(e^{i\theta_1},e^{i\theta_2},e^{i\theta_3})$ if and only if $\theta_1<
\theta_2<\theta_3$. The subgroup $\UU$ of torsion elements of $\KK$, 
that is, the roots of $1$ in the field $\CC$ of all complex numbers, is a cyclically ordered group. 
On the other side, any linearly ordered group $(G,<)$ can be 
equipped with a cyclic order: we set $R(a,b,c)$ if and only if either 
$a < b < c$ or $b < c < a$ or $c < a < b$. 
In this case we say that $(G,R)$ is a {\it linearly cyclically ordered group} (or sometimes 
{\it linear cyclically ordered goup}). \\ 
\indent If $g_0$ is a positive, central and cofinal element of a linearly ordered group $(G,<)$, 
then the quotient group $G/\langle g_0\rangle$ (where $\langle g_0\rangle$ denotes the 
subgroup generated by $g_0$) can be cyclically ordered by setting $R(a\langle g_0\rangle, 
b\langle g_0\rangle,c\langle g_0\rangle)$ if and only if there exist $a'\in a\langle g_0\rangle$, 
$b'\in b\langle g_0\rangle$ and $c'\in c\langle g_0\rangle$ such that $1_G\leq a'<g_0$, 
$1_G\leq b'<g_0$, $1_G\leq c'<g_0$ and either $a'<b'<c'$ or $b'<c'<a'$ or $c'<a'<b'$. This 
cyclically ordered group is called the {\it wound-round} associated to $G$ and $g_0$. 
The theorem of Rieger (\cite{Fu 63}, 
IV, 6, th. 21) states that every cyclically ordered group can be obtained in this way, that is 
for every cyclically ordered group $(G,R)$, there exists a linearly ordered group $\mbox{uw}(G)$ 
and an element $z_G$ which is positive, central and cofinal in $\mbox{uw}(G)$, such that 
$(G,R)$ is isomorphic to the cyclically ordered group $\mbox{uw}(G)/\langle z_G \rangle$. 
$\mbox{uw}(G)$ is the set $\ZZ \times G$ together with the order $(n,c)<(n',c')$ if and only if 
either $n<n'$, or $n=n'$ and $R(1_G,c,c')$, or $n=n'$ and $c=1_G$. We set $z_G=(1,1_G)$, 
and the group law is given by:
$$(m,a)\cdot (n,b)= \left\{ 
\begin{array}{rcl}
(m+n,b)&if&a=1_G\\
(m+n,a)&if&b=1_G\\
(m+n,ab)&if&R(1_G,a,ab)\\
(m+n+1,ab)&if&R(1_G,ab,a)\\
(m+n+1,1_G)&if&ab=1_G\neq a
\end{array}
\right.
$$
$\mbox{uw}(G)$ is called the {\it unwound} of $(G,R)$. 
For example, the unwound of $\KK$ is 
isomorphic to $(\RR,+)$. \\ 
\indent If $C$ is the greatest convex subgroup of the unwound $\mbox{uw}(G)$ of a cyclically 
ordered group $(G,R)$, then $z_G^{-1}<C<z_G$, hence $C$ embeds trivially in $G\simeq 
\mbox{uw}(G)/\langle z_G\rangle$, we denote by $l(G)$ its image and we call it the 
{\it linear part} of $G$. The restriction of $R$ to $l(G)$ is a linearly cyclic order, where the 
order on $l(G)$ is induced by the order on $\mbox{uw}(G)$. \\ 
\indent We will denote by 
$K(G)$ the quotient group $G/l(G)$. The class modulo $l(G)$ of an element $a\in G$ will be 
denoted by $\bar{a}$. $K(G)$ is equipped with a structure of a cyclically ordered group, 
which is induced by the cyclic order on $G$, because if 
$a,\; b,\; c$ in $G$ belong to pairwise distinct classes, and $d\in l(G)$, 
then $R(a,b,c)\Leftrightarrow R(ad,b,c)\& R(da,b,c)$. 
We know that $K(G)$ embeds into the cyclically ordered group $\KK$. \\ 
\indent We let 
$R(x_1,\dots, x_n)$, where $x_1,\dots,\; x_n$ belong to $G$, stand for: 
$$\bigcup_{1\leq i<j<k\leq n} R(x_i,x_j,x_k)\mbox{, which is equivalent to:  }
\bigcup_{1\leq i\leq n-2} R(x_i,x_{i+1},x_n)$$ 

In the following, $(G,R)$ is an abelian cyclically ordered group. 
Since $G$ is abelian, we will use the notation $+$, except in the case of 
subgroups of $\KK$, because it is customary to use the notation $\cdot$. \\ 
\indent 
One can prove the following property by induction. 
\begin{num}{\bf Lemma.}\label{lm0}
Let $g_1,\dots, g_n$ be elements of $G$ and $\theta_1,\dots,\theta_n$ be the elements 
of the interval $[0,2\pi [$ of $\RR$ such that, for $1\leq j\leq n$, $\bar{g}_j=e^{i\theta_j}$. 
Let $k$ be the integer part of $(\theta_1+\cdots + \theta_n)/(2\pi)$ and 
$\theta= (\theta_1+\cdots + \theta_n)-2k\pi$. There exists $g$ in $G$ such that  
$\bar{g}=e^{i\theta}$ and  we have: $(0,g_1)+ \cdots + 
(0,g_n)=(k,g)$ within $\mbox{uw}(G)$. 
\end{num}
\indent We will use this property in order in order to obtain a representation of 
$n(a,x)\in\mbox{uw}(G)$ as a $2$-uple from $\ZZ\times G$. \\ 
\indent Recall that $G$ (a priori non necessarily abelian) is said to be 
{\it c-ar\-chi\-me\-dean} if $l(G)=\{0\}$. 
\'Swierczkowski (\cite{Sw 59}) proved that $G$ is c-archimedean 
if and only if for every $x$ and $y$ in $G$ there exists an integer $n >0$ such that  
$R(0,nx,y)$ does not hold. Note that this is not a first order property. If $G$ is 
c-archimedean, then it embeds into $\KK$, hence it is abelian. 
A {\it c-convex subgroup} of $G$ is a 
subgroup $H$ which satisfies for every $a$ and $b$ in $G$: $(b\in H \Rightarrow b\neq -b)\&
(b\in H\, \&\, R(0,b,-b)\,
\& \,R(0,a,b))\Rightarrow a\in H$. There exists a canonical one-to-one correspondence  
between the set of proper c-convex subgroups of $G$ 
and the set of convex subgroups of $l(G)$ which are different 
from $\{0\}$ (note that $l(G)$ is a c-convex subgroup). 
We see that $G$ is c-archimedean if and only if it doesn't contain any proper 
c-convex subgroup. 
We will say that $G$ is {\it c-$n$-divisible} if it is $n$-divisible and it contains 
a subgroup which is isomorphic to the subgroup of $(\UU,\cdot)$ generated by a 
primitive $n$-th root of $1$. We see that this is a first order property. 
Now, $G$ contains 
a subgroup which is isomorphic to the subgroup of $(\UU,\cdot)$ generated by a 
primitive $n$-th root of $1$ if and only if $z_G$ is $n$-divisible in the unwound 
$\mbox{uw}(G)$. It follows that if $\mbox{uw}(G)$ is $n$-divisible then $G$ is 
c-$n$-divisible. Conversely, let $z\in \mbox{uw}(G)$ be such that $z_G=nz$, then, 
since $G$ is $n$-divisible, for every $x\in \mbox{uw}(G)$ there exists $y\in \mbox{uw}(G)$ 
such that $x-ny$ belongs to the subgroup generated by $z_G$, that is, there exists an integer 
$k$ such that $x-ny=kz_G=nkz$, so $x=n(y+kz)$. Consequently: $G$ is 
c-$n$-divisible if and only if $\mbox{uw}(G)$ is $n$-divisible. 
We will say that $G$ is {\it c-divisible} if it is divisible and it  
contains a subgroup which is isomorphic to $(\UU,\cdot)$. 
\section{c-regular groups.}\label{section2}
\begin{defi}{\bf Definitions.} Let $n\geq 2$ be an integer. We say that 
$(G,R)$ is {\it $n$-regular} if for every $x_1,\dots,x_n$ in $G$ which satisfies  
$R(0,x_1,\dots,x_n,-x_n)$, there exists $x \in G$ such that 
$(R(x_1,nx,x_n)$ or $nx=x_1$, 
or $nx=x_n)$ and $R(0,x,\dots, (n-1)x,x_n)$.  
We say that $(G,R)$ is {\it c-$n$-regular} 
if for every $x_1,\dots,x_n$ in $G$ satisfying $R(0,x_1,\dots,x_n)$, 
there exists $x \in G$ such that 
$(nx=x_1$, or $nx=x_n$, or $R(x_1,nx,x_n))$ and $R(0,x,\dots, (n-1)x,x_n)$. \\ 
We say that $(G,R)$ is {\it regular} if it is $n$-regular 
for every $n\geq 2$, and we say that $(G,R)$ is {\it c-regular} if it is c-$n$-regular 
for every $n\geq 2$. 
\end{defi}
\indent Note that these are first order properties, and clearly if 
$(G,R)$ is c-$n$-regular, then it is $n$-regular. Now from Lemma \ref{lm4b} it 
follows that 
a nontrivial linearly cyclically ordered group is not c-$n$-regular. 
One can prove that if $(G,R)$ is a linearly cyclically ordered group, then 
it is a $n$-regular cyclically ordered group if and only if it is a $n$-regular linearly ordered group, 
that the subgroups of $\KK$ are c-regular. \\ 
\indent The proofs of the three following lemmas are left to the reader. 
\begin{num}{\bf Lemma.}\label{lmap2} Assume that $G$ is a nontrivial subgroup of $\KK$ and 
let $n\geq 2$. Then $G$ is $n$-divisible and contains a primitive $n$-th root of 
$1$ if and only if for every $\theta\in [0,\pi]$, if $e^{i\theta} \in G$, 
then $e^{i\theta/n}\in G$. 
\end{num}
\indent Note that the c-archimedean groups (i.e. the subgroups of $\KK$) are c-archimedean 
but they need not satisfy the conditions of Lemma \ref{lmap2}. 
\begin{num}{\bf Lemma.}\label{lm3} Let $n\geq 2$, assume that $G$ is $n$-regular and is 
not c-archimedean. For every $\theta$ in $[0,\pi]$, 
if $e^{i\theta} \in K(G)$, then $e^{i\theta/n}\in K(G)$.
\end{num}
\begin{num}{\bf Lemma.}\label{lm4} 
Let $n\geq 2$ such that $G$ is $n$-regular, 
then $l(G)$ is $n$-regular, where the compatible total order on $l(G)$ is derived from 
the definition of the positive cone of $G$. In particular, if $G$ is regular, 
then $l(G)$ is regular. 
\end{num}
\begin{num}{\bf Lemma.}\label{lm4b} 
Let $n\geq 2$, assume that $G$ is c-$n$-regular and is not c-ar\-chi\-me\-dean, 
then $K(G)$ contains a primitive $n$-th root of $1$. 
\end{num}
\begin{dem}{\bf Proof.} 
Assume that $h_1<\cdots<h_n<0$ are elements of $l(G)$ such that $2h_n<h_1$. 
In the case where $l(G)$ is dense, clearly one can find such a sequence, and if 
$l(G)$ is discrete with first positive element $\epsilon$, 
let $h_n \in l(G)$ such that $h_n<-n\epsilon$, and set  
$h_1=h_n-n\epsilon+\epsilon, \dots, h_{n-1}=h_n-\epsilon$. Since $G$ is c-$n$-regular, 
there exists $g\in G$ such that $ng=h_1$, or $ng=h_n$, or $R(h_1,ng,h_n)$, and 
$R(0,g,\dots,(n-1)g,h_n)$. 
Note that $\bar{g}$ is a $n$-th root of $1$ in $K(G)$, because $\overline{ng}=\bar{0}$. 
Now $l(G)$ is c-convex, so 
$ng\in l(G)$, and $h_1\leq ng \leq h_n<0$. Let $k$, $1\leq k\leq n$, be such that  
$kg\in l(G)$, and $d$ be the gcd of $k$ and $n$. Let 
$a$ and $b$ be two integers such that 
$ak+bn=1$, then $dg=akg+bng\in l(G)$, hence, if $k$ is the lowest positive integer 
such that $kg\in l(G)$, then $k$ divides $n$. If $k<n$, since $R(0,kg,h_n)$ holds, 
we have $kg<h_n<0$.  Furthermore $2k\leq n$, hence 
$ng\leq 2kg <2h_n<h_1$: a contradiction. Consequently, $k=n$. 
Since for $1\leq k\leq n-1$ $\overline{kg}\neq \bar{0}$, $\bar{g}$ is a primitive $n$-th root. \end{dem}
\begin{num}{\bf Proposition.}\label{prop1} Let $n\geq 2$, assume that $l(G)$ is 
$n$-regular. If $K(G)$ is $n$-divisible and contains a primitive $n$-th root
of $1$, then $G$ is c-$n$-regular. 
\end{num}
\begin{dem}{\bf Proof.} 
$K(G)$ is nontrivial because it contains a primitive $n$-th root 
of $1$; since it is $n$-divisible, 
it is infinite, hence it is a dense subgroup of $\KK$. \\ 
\indent 
Let $g_1$ and $g_2$ in $G$, assume: $R(0,g_1,g_2)$ and $\bar{g_1}\neq \bar{g_2}$. 
There exist $\theta_1<\theta_2$ in $[0,2\pi]$ such that $\bar{g_1}=e^{i\theta_1}$ and 
$\bar{g_2}=e^{i\theta_2}$. Since $K(G)$ is dense, there exist $g$ in $G$ and $\theta$ 
in $]\theta_1/n,\theta_2/n[$ such that $\bar{g}=e^{i\theta}$. Then 
$\overline{ng}=e^{in\theta}$, with $0\leq \theta_1<n\theta<\theta_2$, and 
$0<\theta<\cdots < (n-1)\theta <\theta_2$, hence $R(g_1,ng,g_2)$ and 
$R(0,g,\dots, (n-1)g,g_2)$. \\ 
\indent Now, let $g_1,\dots, g_n$ in $G$ be such that  
$R(0,g_1,\dots, g_n)$. If $\bar{g_1}\neq \bar{g_n}$, then by what has been proved above, 
there exists $g\in G$ such that $R(g_1,ng,g_n)$ and 
$R(0,g,\dots, (n-1)g,g_n)$. If $\bar{g_1}= \bar{g_n}=\bar{0}$ and $g_1\in l(G)_+$, 
$g_n\in -l(G)_+$, we set $\bar{g_1}=e^{i0}$ and $\bar{g_n}=e^{2i\pi}$, and 
we conclude in the same way. We assume now that $\bar{g_1}=\cdots=\bar{g_n}$, 
we set $\bar{g_1}=e^{i\theta}$, with $0\leq \theta \leq 2\pi$, 
and if $\theta=0$ we assume that either each of 
$g_1,\dots,g_n$ belongs to $l(G)_+$ or each of $g_1,\dots,g_n$ belongs to $-l(G)_+$. 
In the first case, we set $\theta	=0$, in the second case, we set 
$\theta=2\pi$. 
Note that in any case $l(G)$ is nontrivial. \\ 
\indent If $\theta=0$, i.e. both of 
$g_1,\dots,g_n$ belong to $l(G)_+$, then we have $0<g_1<\cdots 
<g_n$. Since $l(G)$ is $n$-regular, there exists $g\in l(G)$ such that $g_1\leq ng\leq g_n$, and 
we have also $0<g<\cdots< (n-1)g<g_n$, hence $(R(g_1,ng,g_n)$ or $ng=g_1$ or 
$ng=g_n)$ and $R(0,g,\dots,(n-1)g,g_n)$.  \\
\indent Now assume that $0<\theta \leq 2\pi$. 
Since $K(G)$ is $n$-divisible and it contains a primitive $n$-th root of $1$, 
there exists an 
element $g$ of $G$ such that $\bar{g}=e^{i\theta/n}$. 
Since $g_1-ng, \dots, g_n-ng$ belong to 
$l(G)$, it follows from the relation $R(0,g_1,\dots, g_n)$ that $g_1-ng<\cdots <
g_n-ng$ in $l(G)$. After replacing $g$ by $ng-g_1$, if necessary, we can assume 
that $0<g_1-ng<\cdots<g_n-ng$. Since $l(G)$ is $n$-regular, 
there exists $h\in l(G)$ such that $g_1-ng\leq nh\leq g_n-ng$. It follows: 
$R(g_1-ng,nh,g_n-ng)$, or $nh=g_1-ng$, or $nh=g_n-ng$, 
then by compatibility $R(g_1,n(g+h),g_n)$ or $n(g+h)=g_1$ or $n(g+h)=g_n$. Furthermore, 
since $0<\theta/n<\cdots <(n-1)\theta/n<\theta$, we have 
$R(0,g+h,\dots,(n-1)(g+h),g_n)$. \\ 
\indent This proves that $G$ is c-$n$-regular. 
\end{dem}
\begin{num}{\bf Theorem.}\label{prop2} Assume that $G$ is not c-archimedean 
and let $n\geq 2$. The following assertions are equivalent. \\ 
1) $G$ is $n$-regular and $K(G)\neq \{0\}$. \\
2) $l(G)$ is $n$-regular,  $K(G)\neq \{0\}$ and for every 
$\theta \in [0,\pi]$ such that $e^{i\theta}\in K(G)$ we have: $e^{i\theta/n}\in K(G)$. \\ 
3) $G$ is c-$n$-regular. \\
4) $l(G)$ is $n$-regular, $K(G)$ is $n$-divisible and contains a primitive 
$n$-th root of $1$. 
\end{num}
\begin{dem}{\bf Proof.} 1) $\Rightarrow$ 2). 
If $G$ is $n$-regular, then by Lemma 
\ref{lm3}, for every $\theta \in [0,\pi]$ such that $e^{i\theta}\in K(G)$ 
we have $e^{i\theta/n}\in K(G)$, furthermore $K(G)$ is $n$-divisible, and by Lemma \ref{lm4}, $l(G)$ is $n$-regular. 
4) $\Rightarrow$ 3)  follows from Proposition \ref{prop1}. 3) $\Rightarrow$ 1) follows from  
the definition, and from Lemma \ref{lmap2} we deduce 2) $\Leftrightarrow$ 4). 
\end{dem}
\begin{num}{\bf Theorem.}\label{prop2b} Let $n\geq 2$. $G$ is $n$-regular if and only if 
one of the two following conditions is satisfied. \\ 
1) $G$ is linearly cyclically ordered, and is a $n$-regular linearly ordered group. \\ 
2) $G$ is not linearly cyclically ordered, and is c-$n$-regular. \\ 
In other words, the class of $n$-regular cyclically ordered groups is the 
union of the class $n$-regular linearly cyclically ordered groups and of the 
class of c-$n$-regular cyclically ordered groups. 
\end{num}
\begin{dem}{\bf Proof.} Straightforward. 
\end{dem} 
\begin{num}{\bf Lemma.}\label{lm5} 
Let $n\geq 2$, assume that there exist $g\in G$ and $x\in 
\mbox{uw}(G)$ such that $\bar{g}=1$ and $nx=(1,g)$, then $e^{2i\pi/n}\in K(G)$. 
It follows that if $l(G)$ is nontrivial and $\mbox{uw}(G)$ is $n$-regular, 
then $e^{2i\pi/n} \in K(G)$. 
\end{num}
\begin{dem}{\bf Proof.} Since $(1,g)>0$, by compatibility we have: $0<x<(1,g)$. By the 
definition of the order on $\mbox{uw}(G)$, there exists $h\in G$ such that either $x=(0,h)$ or 
$x=(1,h)$. Now $n(1,h)\geq (n,nh)\geq (n,0)>(1,g)$, hence $x=(0,h)$. Set  
$\bar{h}=e^{i\theta}$, with $0\leq \theta <2\pi$. 
According to Lemma \ref{lm0}, we have $n\theta =2\pi$, which proves that  
$e^{2i\pi/n}\in K(G)$. Now assume that $l(G)$ is 
nontrivial and that $\mbox{uw}(G)$ is $n$-regular. 
Let $g_1,\dots, g_n$ in $l(G)$, assume that $0<g_1<\cdots <g_n$, hence 
$(0,0)<(1,g_1)<\cdots <(1,g_n)$. There exists $y\in \mbox{uw}(G)$ such that 
$(1,g_1)\leq ny\leq (1,g_n)$. Then $ny$ admits a representation as $ny=(1,g)$ with $g\in G$. 
By the definition of the order on $\mbox{uw}(G)$, $R(g_1,g,g_n)$ holds, 
and by c-convexity of $l(G)$, we have $g\in l(G)$. 
Finally, from what we just proved, $e^{2i\pi/n}\in K(G)$. 
\end{dem} 
\begin{num}{\bf Lemma.}\label{lm6} Let $n\geq 2$, if $\mbox{uw}(G)$ is 
$n$-regular, then $G$ is $n$-regular. 
\end{num}
\begin{dem}{\bf Proof.} 
Let $g_1,\dots,g_n$ in $G$, assume that $R(0,g_1,\dots,g_n)$ holds. Since $\mbox{uw}(G)$ 
is $n$-regular, and $(0,0)<(0,g_1)<\cdots < (0,g_n)$, there exists $x$ in $\mbox{uw}(G)$ 
such that $(0,g_1)\leq nx\leq (0,g_n)$. In the same way of in the proof of Lemma 
\ref{lm5}, $nx$ admits a representation as $nx=(0,g)$ with $g\in G$. 
It follows that $x$ admits a representation as  
$x=(0,h)$ with $nh=g$. We have $(0,g_1)\leq (0,nh)\leq (0,g_n)$, hence $R(g_1,nh,g_n)$ 
or $nh=g_1$, or $nh=g_n$, which proves that $G$ is $n$-regular. 
\end{dem} 
\begin{num}{\bf Proposition.}\label{prop3} Let $n\geq 2$, assume that 
$l(G)$ is $n$-regular, and that $K(G)$ is $n$-divisible and contains a primitive 
$n$-th root of $1$. 
Then $\mbox{uw}(G)$ is $n$-regular. 
\end{num}
\begin{dem}{\bf Proof.} By assumption, $K(G)$ is nontrivial, hence $G$ is not equal to $l(G)$. 
Hence if $G$ is c-archimedean, then $\mbox{uw}(G)$ 
is archimedean, hence it is regular. If $G$ is not c-archimedean, then $l(G)$ is 
nontrivial, hence $G$ is infinite. $K(G)$ is infinite because it is $n$-divisible. 
Let $(0,0)<(m_1,g_1)<\cdots<(m_n,g_n)$ in $\mbox{uw}(G)$. 
Let $q$ and $r$ be the integers such that $0\leq r<n$ and 
$m_1=nq+r$. If we can find $x\in \mbox{uw}(G)$ such that 
$(r,g_1)\leq nx\leq (m_n-nq,g_n)$, then 
we have $(m_1,g_1)\leq n(x+(q,0))\leq (m_n,g_n)$, which proves that we can assume that  
$0\leq m_1<n$. Since $K(G)$ contains a primitive $n$-th root of $1$, we have: 
$e^{2i\pi/n}\in K(G)$. 
Let $g\in G$ such that $\bar{g}=e^{2i\pi/n}$, and set $g'=m_1g$, then 
$\bar{g'}=e^{2im_1\pi/n}$, with $0\leq 2m_1\pi/n<2\pi$, and $n(2m_1\pi/n)=2m_1\pi$. 
By Lemma \ref{lm0}, there exists $h\in l(G)$ such that $n(0,g')=
(m_1,h)$. If $\bar{g_1}=e^{i\theta_1}$, since $K(G)$ is divisible, one can find 
an element whose class is $e^{i\theta_1/n}$ in $K(G)$. We add such an element to $g'$, 
then, if necessary, we add an element of $l(G)$ in order to get $g''$ such that 
$n(0,g'')=(m_1,ng'')$ and $R(0,ng'',g_1)$ holds. 
Assume that $m_1=\cdots =m_n$ and $\bar{g_1}=\cdots=\bar{g_n}$. 
By substracting $n(0,g'')$ to each term, we get a sequence 
$(0,0)<(0,h_1)<\cdots<(0,h_n)$, where the $h_k$'s belong to $l(G)$. By the  
definition of the order on $\mbox{uw}(G)$, we have either $0<h_1<\cdots<h_n$, 
or $h_1<\cdots <h_n<0$, or $h_n<0<h_1$. 
In the first case, there exists $h$ in 
$l(G)$ such that $h_1\leq nh\leq h_n$. Hence $R(h_1,nh,h_n)$, or $nh=h_1$, or $nh=h_n$. 
In the second case we consider the sequence 
$0<-h_n<\cdots<-h_1$, there exists $h\in l(G)$ such that $-h_n\leq -nh\leq -h_1$, then 
$h_1\leq nh \leq h_n$, and we conclude in the same way as above. In the third case, 
we set for example $h=h_1$, then $h_n<0<h_1<nh$, which implies  
$R(h_1,nh,h_n)$. In any case, we have proved that 
there exists $h\in l(G)$ such that $h_1\leq nh\leq h_n$. Then 
$(0,h_1)\leq (0,nh)=n(0,h)\leq (0,h_n)$. 
Consequently $(m_1,g_1)\leq n((0,h)+(0,g''))\leq (m_n,g_n)$. 
If the $m_k$'s are not equal, note that we have $m_1\leq 
\cdots \leq m_n$, and in this case $m_1<m_n$. Let $0<h_2<\cdots<h_n$ in 
$l(G)$, then $0<(m_1,g_1)<(m_1,g_1+h_2)<\cdots<(m_1,g_1+h_n)$. It follows that 
there exists $h\in G$ such that $(m_1,g_1)\leq n(0,h)\leq (m_1,g_1+h_n)$, and 
since $(m_1,g_1+h_n)<
(m_n,g_n)$, $h$ is the required element. If all of the $m_k$'s are equal, but the 
$g_k$'s are not, then $R(0,g_1,g_n)$ holds, take the element $(0,h)$ that we just exhibit 
above. 
\end{dem}
\begin{num}{\bf Theorem.}\label{prop4} Let $n\geq 2$. The following conditions 
are equivalent. \\ 
(1) $\mbox{uw}(G)$ is $n$-regular \\ 
(2) $G$ is c-$n$-regular\\
(3) either $G$ is c-archimedean, \\ 
or 
$l(G)$ is $n$-regular, $K(G)$ is $n$-divisible and it contains a primitive 
$n$-th root of $1$ \\ 
(4) the quotient of $G$ by every proper c-convex subgroup is c-$n$-divisible.
\end{num}
\begin{dem}{\bf Proof.} If $G$ is c-archimedean, then it is either dense or finite, hence 
it is c-regular. Furthermore, 
$\mbox{uw}(G)$ is archimedean, hence it is regular. Now, assume that  
$G$ is not c-archimedean. The equivalence (2) $\Leftrightarrow$ (3) has been proved 
in Theorem \ref{prop2}. 
If $l(G)$ is $n$-regular, $K(G)$ is $n$-divisible and contains 
a primitive $n$-th root of $1$, then 
by Proposition \ref{prop3}, $\mbox{uw}(G)$ is $n$-regular. 
If $\mbox{uw}(G)$ is $n$-regular, since it contains a convex subgroup 
which is isomorphic to $l(G)$, $l(G)$ is $n$-regular. 
According to Lemma \ref{lm6}, $G$ is $n$-regular, 
and by Lemma \ref{lm3}, $K(G)$ is $n$ divisible. 
Finally, by Lemma \ref{lm5}, $e^{2i\pi/n}\in K(G)$. \\ 
\indent Let $H$ be a proper c-convex subgroup (that is a convex subgroup of $l(G)$ 
distinct from $\{0\}$). In the construction of the unwound, 
we note that $\mbox{uw}(G/H)=\mbox{uw}(G)/H$. By means of 
(1) $\Leftrightarrow$ (2) we get: \\ 
$G/H$ is c-$n$-divisible, for every c-convex subgroup $H$ of $G$,  if and only if 
$\mbox{uw}(G)/H$ is $n$-divisible, for every c-convex subgroup $H$ of $G$, 
if and only if $\mbox{uw}(G)$ is $n$-regular, 
if and only if $G$ is c-$n$-regular. 
\end{dem}
\section{Elementarily equivalent c-regular cyclically ordered groups.}\label{section3} 
\indent As being dense and being c-regular are first order properties, 
every cyclically ordered group which is elementarily equivalent to 
an  infinite subgroup of $\KK$ is dense 
and c-regular. Every ultraproduct of finite cyclically ordered groups 
is discrete and c-regular, since each factor satisfies both of these first order properties. \\ 
\indent If $G$ is discrete and infinite then $l(G)$ is nontrivial, because 
every infinite and c-archimedean cyclically ordered group is dense. 
Furthermore, if $G$ is c-regular, then by Lemma \ref{lm3}, $K(G)$ is divisible, hence, 
it is infinite. 
\subsection{Preliminaries.}\label{subsec31}
If $A$ is an abelian group and $p$ is a prime, we define the $p$-th 
{\it prime invariant of Zakon} of $A$, denoted by $[p]A$, to be the maximum number 
of $p$-incongruent elements in $A$. In the infinite case, we set $[p]A=\infty$, 
without distinguishing between infinities of different cardinalities. 
\begin{num}{\bf Lemma.}\label{lm33} 
Let $H$ be a subgroup of $\QQ$, and, for every prime $p$, 
let $m_p\in \NN$ be such that $[p]H\leq p^{m_p}$. Then there exists a countable subgroup $M$ 
of $(\RR,+)$ which contains $H$, such that $H$ is pure in $M$ and, for every prime 
$p$, $[p]M=p^{m_p}$. 
\end{num}
\begin{dem}{\bf Proof.} 
Denote by $p_1<p_2<\cdots <p_n<\cdots$ the increasing sequence of all primes, 
and, for every 
$n\in \NN^*$, denote by $m_n'$ the integer which satisfies $p_n^{m_n'}=P_n^{m_{p_n}}/[p_n]H$. 
We pick a transcendent number $e$, and, for $n$ and $j$ in $\NN^*$ such that  
$1\leq j \leq m_n'$, we pick pairwise distinct elements $k_{n_j}$ in  
$\NN^*$. Denote by $M_0$ the direct sum of the $\ZZ_{(p_n)}e^{k_{n_j}}$'s, where 
$n\in \NN^*$ and $1\leq j \leq m_n'$, and $\ZZ_{(p_n)}$ denotes the localization of $\ZZ$ at the 
prime ideal $(p_n)$. By \cite{Za 61}, for every $n\in\NN^*$ we have: $[p_n]M_0=p^{m_n'}$. 
Since $\QQ\cap M_0=\{0\}$, we have: $\QQ+M_0=\QQ\oplus M_0$. We set: 
$M=H\oplus M_0$, then $M$ contains 
$H$, $H$ is pure in $M$, $M$ is countable, and for every prime $p$ we have: 
$[p]M=([p]H)([p]M_0)=p^{m_p}$. 
\end{dem}
\begin{num}{\bf Proposition.}\label{propteta1} 
Let $A$, $B$, $T$ be linearly ordered abelian groups, where $B$ is a convex 
subgroup of $T$. Then there exists an exact sequence of ordered groups 
$0\rightarrow B \rightarrow T \rightarrow A \rightarrow 0$ if and only if  
there exists a mapping $\theta$ from $A\times A$ to $B$ 
which satisfies: \\
(*) $\forall a\in A,\; \theta(a,0)=0$, \\ 
(**) $\forall (a_1,a_2) \in A \times A,\; \theta(a_1,a_2)=\theta(a_2,a_1)$, \\ 
(***) $\forall (a_1,a_2,a_3)\in A\times A\times A,\; \theta(a_1,a_2)+
\theta(a_1+a_2,a_3)=\theta(a_1,a_2+a_3)+\theta(a_2,a_3)$ (in other words, $\theta$ 
is a $2$-cocycle), \\ 
such that $T$ is isomorphic to $A\times B$ lexicographically ordered and equipped 
with the operation $+_{\theta}$ defined by: 
$$\forall(a_1,b_1,a_2,b_2)\in A\times B\times A\times B,\; 
(a_1,b_1)+_{\theta}(a_2,b_2)=(a_1+a_2,b_1+b_2+\theta(a_1,a_2)).$$
We will say that this exact sequence is an extension of $B$ by $A$, we set 
$T=A\overrightarrow{\times}_{\theta}B$, and we will omit $\theta$ when 
$\theta \equiv 0$. 
\end{num}
\begin{dem}{\bf Proof.} See \cite{Ja 54}. We recall the definition of the mapping $\theta$. 
For every $a\in A$, we pick some $r(a)$ 
in $T$  whose image is $a$ ($a$ can be viewed as a class modulo $B$), 
by taking $r(0)=0$, and we set, for every $a_1$, $a_2$ in $A$, $\theta(a_1,a_2)=r(a_1+a_2)-r(a_1)-r(a_2)$. 
The condition (*) is satisfied because we set $0=r(0)$, (**) follows from commutativity, and 
(***) follows from the associativity of addition. We get the isomorphism by defining the image 
of $(a,b)$, where $a\in A,\; b\in B$, to be the element $r(a)+b$ of $T$. 
\end{dem}
\begin{num}{\bf Proposition.}\label{propteta1b} \\ 
1) Let $A\overrightarrow{\times}_{\theta_1}B_1 $ and 
$A\overrightarrow{\times}_{\theta_2}B_2$ be two extensions of 
linearly ordered abelian groups, $B=B_1\overrightarrow{\times} B_2$, 
$\theta=(\theta_1,\theta_2)$ from $A\times A$ to $B$, then 
$A\overrightarrow{\times}_{\theta}B$ is an extension of 
linearly ordered abelian groups. 
Furthermore $A\overrightarrow{\times}_{\theta_1}B_1$ and 
$A\overrightarrow{\times}_{\theta_2}B_2$ embed canonically into cofinal 
subgroups of $A\overrightarrow{\times}_{\theta}B$. \\
2) Let $A_1\overrightarrow{\times}_{\theta_1}B$ and 
$A_2\overrightarrow{\times}_{\theta_2}B$ be two extensions of 
linearly ordered abelian groups, $A=A_1\overrightarrow{\times} A_2$, 
$\theta=\theta_1+\theta_2$ from $A\times A$ to $B$, then 
$A\overrightarrow{\times}_{\theta}B$ is an extension of 
linearly ordered abelian groups. 
Furthermore $A_1\overrightarrow{\times}_{\theta_1}B$ and 
$A_2\overrightarrow{\times}_{\theta_2}B$ 
embed canonically into subgroups of
$A\overrightarrow{\times}_{\theta}B$. \\
3) Let $A$ and $B'$ be two linearly ordered abelian groups, 
$B$ be a subgroup of $B'$, and 
$A\overrightarrow{\times}_{\theta}B$ be an extension of linearly ordered abelian groups.
Then there exists an extension of linearly ordered abelian groups  
$A\overrightarrow{\times}_{\theta'}B'$, where $\theta'$ extends $\theta$, 
$A\overrightarrow{\times}_{\theta}B$ is a subgroup 
of $A\overrightarrow{\times}_{\theta'}B'$, $\{0\}\overrightarrow{\times}_{\theta'}B'$ 
is a convex subgroup of 
$A\overrightarrow{\times}_{\theta'}B'$, and the quotient group is isomorphic to $A$. 
\end{num}
\begin{dem}{\bf Proof.} The proof is left to the reader. 
\end{dem}
\indent Note that Proposition \ref{propteta1} remains true by assuming that  
$A$ and $T$ are cyclically ordered, and $B$ is linearly ordered. Then analogues 
of 1) and 3) of Proposition \ref{propteta1b} hold. 
\begin{num}{\bf Proposition.}\label{propteta2} 
Let $A'$ and $B$ be two linearly ordered abelian groups, 
$A$ be a subgroup of $A'$, and 
$A\overrightarrow{\times}_{\theta}B$ be an extension of linearly ordered abelian groups.
Then $\theta$ extends by induction to $A' \times A'$ by means of the following 
properties. Let $a' \in A'\backslash A$. \\ 
1) If $a'$ is rationally independent from $A$, we can extend $\theta$ to $\ZZ a'+A$ ($=\ZZ a'\oplus A$) in the following way: 
$$\forall (n_1,a_1,n_2,a_2)\in \ZZ\times A\times \ZZ\times A,\; 
\theta (n_1a'+a_1,n_2a'+a_2)=\theta (a_1,a_2).$$
2) If there exists a prime $p$ such that $pa'\in A$, then for every 
$b_0\in B$, $\theta$ extends to the subgroup generated by $a'$ and $A$ in the 
following way. The elements of the group generated by $a'$ and $A$ are the $na'+a$, 
where $0\leq n \leq p-1$ and $a\in A$. Denote by $\mbox{int}(x)$ the integer part of a 
rational number $x$. The extension is defined by:  
$$\forall (n_1,a_1,n_2,a_2)\in \{0,\dots , p-1\}\times A\times \{0,\dots ,p-1\}\times A,$$ 
$$\theta(n_1a'+a_1,n_2a'+a_2)=\theta(a_1,a_2)+\theta	(\mbox{int}(\frac{n_1+n_2}{p})
(pa'),a_1+a_2)+\mbox{int}(\frac{n_1+n_2}{p})b_0.$$ 
So we get $[p]B$ non isomorphic extensions of $(A\times B,+_{\theta})$, 
each one depending on the class of $b_0$ modulo $pB$. 
$A'\times B$ is a subgroup of the divisible closure of 
$A\overrightarrow{\times}_{\theta}B$. 
If $b_0=0$, then $p(a',0)=(pa',0)$. 
\end{num}
\begin{dem}{\bf Proof.} In case 1) as well as in case 2), one easily checks that 
the mapping $\theta$ that we define satisfies (*) and (**). In case 1), one can 
easily verify that (***) holds, it remains to prove (***) in case 2). 
In order to simplify the notations, we set:
$$ \epsilon_{12}=\mbox{int}(\frac{n_1+n_2}{p}), \;
\epsilon_{23}=\mbox{int}(\frac{n_2+n_3}{p}) \mbox{ and }
\epsilon_{123}=\mbox{int}(\frac{n_1+n_2+n_3}{p}).$$
So let $n_1$, $n_2$, $n_3$ be in $\{0,\dots ,p-1\}$ and $a_1$, $a_2$, $a_3$ 
be in $A$. In order to compute $\theta(n_1a'+a_1+n_2a'+a_2,n_3a'+a_3)$, 
we start from the representation of $n_1a'+a_1+n_2a'+a_2$ as 
$$(n_1+n_2-\epsilon_{12})(pa')+a_1+a_2+
\epsilon_{12}(pa'),$$ and we get: 
$$\theta(n_1a'+a_1+n_2a'+a_2,n_3a'+a_3) 
= \theta(a_1+a_2+\epsilon_{12}(pa'),a_3)+$$ 
$$\theta(\mbox{int}((n_1+n_2-\epsilon_{12}p+n_3)/p)(pa'),
a_1+a_2+a_3+\epsilon_{12}(pa'))$$ 
$$+ \mbox{int}((n_1+n_2-\epsilon_{12}p+n_3)/p)b_0= 
\theta(a_1+a_2+\epsilon_{12}(pa'),a_3)+$$ 
$$
\theta((\epsilon_{123}-\epsilon_{12})(pa'),
a_1+a_2+a_3+\epsilon_{12}(pa'))+ 
(\epsilon_{123}-\epsilon_{12})b_0.$$ 
In the same way, 
$$\theta(n_1a'+a_,n_2a'+a_2+n_3a'+a_3)=
\theta(a_1,a_2+a_3+\epsilon_{23}(pa'))+$$ 
$$
\theta((\epsilon_{123}-\epsilon_{23})(pa'),
a_1+a_2+a_3+\epsilon_{23}(pa'))+ 
(\epsilon_{123}-\epsilon_{23})b_0.$$ 
So, proving (***) reduces to proving that the following $(E1)$ and $(E2)$ are equal: 
$$(E1)=\theta(a_1+a_2+\epsilon_{12}(pa'),a_3)+
\theta((\epsilon_{123}-\epsilon_{12})(pa'),
a_1+a_2+a_3+\epsilon_{12}(pa'))+$$ 
$$(\epsilon_{123}-\epsilon_{12})b_0+\theta(a_1,a_2)
+\theta(\epsilon_{12}(pa'),a_1+a_2)+\epsilon_{12}b_0$$
$$\mbox{and } (E2)=\theta(a_2,a_3)+\theta(\epsilon_{23}(pa'),a_2+a_3)
+\epsilon_{23}b_0+ \theta(a_1,a_2+a_3+\epsilon_{23}(pa'))+$$ 
$$
\theta((\epsilon_{123}-\epsilon_{23})(pa'),
a_1+a_2+a_3+\epsilon_{23}(pa'))
+ 
(\epsilon_{123}-\epsilon_{23})b_0.$$ 
One can check that we have: 
$$(E1) =
\theta(a_1,a_2)+\theta(a_1+a_2,\epsilon_{12}(pa'))+
\theta(a_1+a_2+\epsilon_{12}(pa'),a_3)$$ 
$$+\theta(a_1+a_2+
\epsilon_{12}(pa')+a_3,(\epsilon_{123}-
\epsilon_{12})(pa'))+\epsilon_{123}b_0.$$  
$$\mbox{and } (E2) =\theta(a_2,a_3)+\theta(a_2+a_3,\epsilon_{23})(pa'))+
\theta(a_2+a_3+\epsilon_{23}(pa'),a_1)$$ 
$$+\theta(a_2+a_3+
\epsilon_{23}(pa')+a_1,(\epsilon_{123}-
\epsilon_{23})(pa'))+\epsilon_{123}b_0.$$  
The operation $+_{\theta}$ being associative, we note that 
$$[[[(a_1,0)+_{\theta}(a_2,0)]+_{\theta}(\epsilon_{12}(pa'),0)]
+_{\theta}(a_3,0)]+_{\theta}
((\epsilon_{123}-\epsilon_{12})(pa'),0)=$$ 
$$
(a_1+a_2+\epsilon_{12}(pa')+a_3+
(\epsilon_{123}-\epsilon_{12})(pa'),(E1)),\mbox{ and that}$$ 
$$[[[(a_2,0)+_{\theta}(a_2,0)]+_{\theta}(\epsilon_{23}(pa'),0)]
+_{\theta}(a_1,0)]+_{\theta}
((\epsilon_{123}-\epsilon_{23})(pa'),0)=$$ 
$$
(a_2+a_3+\epsilon_{23}(pa')+a_1+
(\epsilon_{123}-\epsilon_{23})(pa'),(E2)).$$ 
Now, calculations show that in order to 
prove that $(E1)=(E2)$ it is sufficient to establish: 
$$(\epsilon_{12}(pa'),0)+_{\theta}
(\epsilon_{123}-\epsilon_{12})(pa'),0)= 
(\epsilon_{23}(pa'),0)+_{\theta}
(\epsilon_{123}-\epsilon_{23})(pa'),0).$$ 
The $3$-tuple  
$(\epsilon_{123},\epsilon_{12},
\epsilon_{23})$ can take the values  
$$(0,0,0),\; (2,1,1),\; (1,0,0), \;(1,1,0),\; (1,0,1)\mbox{ or } (1,1,1),$$ 
in any case, calculations prove that the equality holds, the details are left  to the reader. \\ 
\indent If $n_1+n_2<p$, then $(n_1a',0)+_{\theta}(n_2a',0)=((n_1+n_2)a',0)$, and 
$\theta((p-1)a',a')=\theta(0,0)+\theta(pa',0)+b_0=b_0$. Hence 
$((p-1)a',0)+_{\theta}(a',0)=(pa',b_0)$ i.e. $p(a',0)=(pa',b_0)$. In particular, if $b_0=0$, 
then $p(a',0)=(pa',0)$, and 
$(pa',0)$ is divisible by $p$ in the extension we obtained. 
In all cases, $(a',0)$ belongs to the divisible hull of 
$A\overrightarrow{\times}_{\theta}B$. \\ 
\indent Now, the order relation on the quotient group is fixed, hence any morphism 
from $A'\times B$ onto $A'\times B$, equipped with another operation, takes an element 
into an element 
which belongs to the same class modulo $B$. In particular, the image of $(a',0)$ 
is $(a',b_1)$, for some $b_1\in B_1$. Then $p(a',b_1)=p(a',0)+p(0,b_1)=
(pa',b_0)+(0,pb_1)=(pa',b_0+pb_1)$. Consider two operations extending the one of 
$A\times B$, $+_{\theta}$ and $+_{\theta'}$, characterized respectively by 
 elements of $B$, say $b_0$ and $b_0'$. So, if there exists a morphism between these 
two extensions, then $b_0'$ can be expressed as $b_0+pb_1$, for some  $b_1\in B_1$, 
and conversely, which proves that 
there are (at least) the same number of non isomorphic extensions of the structure 
of group as the number of classes modulo $pB$ in $B$. \\ 
\indent If $a'\in A'\backslash A$ satisfies some relation $na'\in A$, 
we consider $n$ as factored into prime powers, 
and we proceed in the same way as in 2). So, $\theta$ extends to every 
element of $A'\backslash A$, and by Zorn's lemma we conclude that 
$\theta$ extends to $A'$. $A'\overrightarrow{\times}_{\theta}B$ is contained in the 
divisible hull of $A\overrightarrow{\times}_{\theta}B$, 
because this holds at each step. 
\end{dem} 
\begin{defi}{\bf Definition.}\label{defteta5} 
Let $T$ be a discrete linearly ordered abelian group, with first positive element 
$1_T$, and which contains a fixed element $z_T$ that we join to the  
language. Since $1_T$ is definable, we can assume that it lies in the language. 
For every prime $p$, $n\in \NN^*$ and $k\in \{0,\dots,p^n-1\}$, 
we define the formula  
$DD_{p^n,k}$: $\exists x,\; p^n x=z_T+k1_T$. 
\end{defi}
\begin{num}{\bf Lemma.}\label{lemteta5b} 
Let $T$ be a discrete linearly ordered abelian group, with first positive element $1_T$, 
and containing a fixed cofinal positive element $z_T$, 
that we join to the language. Assume that $T/\langle 1_T\rangle$  is divisible. \\ 
1) For every prime $p$ and $n\in \NN^*$, there exists exactly one integer $k\in \{0,
\dots , p^n-1\}$ such that $DD_{p^j,k}$ holds in $T$. \\ 
2) Let $p$ be a prime, $n\in \NN^*$, and $k\in \{0,\dots , p^n-1\}$ such that $DD_{p^j,k}$ 
holds within $T$. We express $k$ as a sum $a_0+a_1p+\cdots + 
a_{n-1}p^{n-1}$, where, for $1\leq j \leq n$, $a_j\in \{0,\dots, p-1\}$, then for every 
$j\leq n$, and $k_j=a_0+\cdots+a_{j-1}p^{j-1}$, $DD_{p^j,k_j}$ holds within $T$. 
\end{num}
\begin{dem}{\bf Proof.} 1) Since $T/\langle 1_T\rangle$ is divisible, 
there exists $x \in T$ such that the class of  
$p^nx$ modulo $\langle 1_T\rangle$ is equal to the class of $z_T$. Let $x_1\in T$, 
the class of $x_1$ is equal to the class of $x$ if and only if $x_1-x\in \ZZ\cdot1_T$, 
and if this holds then $p^nx_1-p^nx \in p^n\ZZ\cdot1_T$. So there exists a unique element of 
the class of $x$ (that we still denote by $x$) such that $p^nx-z_T\in \{0,\dots, p^n-1\}$, 
which proves 1). \\
\indent 2) If $p^n x=z_T+k1_T=z_T+(a_0+a_1p+\cdots + a_{n-1}p^{n-1})1_T$, then for 
$j\leq n$ we have $p^j(p^{n-j}x-(a_j+a_{j+1}p+\cdots +a_{n-1}p^{n-j-1})1_T)=
z_T+(a_0+\cdots + 
a_{j-1}p^{j-1})1_T$, the proposition follows. 
\end{dem}
\begin{num}{\bf Proposition.}\label{propteta6} 
Let $T_1$ and $T_2$ be discrete linearly ordered abelian groups such that 
$T_1/\langle1_{T_1} \rangle \simeq \QQ\simeq 
T_2/\langle1_{T_2}\rangle $, and containing an element $z_1$ and $z_2$ respectively, 
which is positive and cofinal 
(which is equivalent to saying that it is not contained in any proper convex subgroup). 
In the language of ordered groups together with a predicate for the distinguished element, the following conditions are equivalent. \\ 
$T_1$ and $T_2$ are isomorphic \\ 
$T_1$ and $T_2$ satisfy the same formulas $DD_{p^j,k}$ \\ 
$T_1 \equiv T_2$.  
\end{num} 
\begin{dem}{\bf Proof.} Clearly, if $T_1$ and $T_2$ are isomorphic, then they are elementarily 
equivalent, and if they are elementarily equivalent then they satisfy the same formulas $DD_{p^j,k}$, because these are first order formulas. 
Now, assume that they satisfy the same formulas $DD_{p^j,k}$. We are going to 
define cocycles $\theta_1$ and 
$\theta_2$ such that $T_1\simeq \QQ \overrightarrow{\times}_{\theta_1} \ZZ$ and 
$T_2\simeq \QQ \overrightarrow{\times}_{\theta_2} \ZZ$, then we will prove: $\theta_1=
\theta_2$. We know that every element of $\QQ$ has a unique representation as a finite sum 
$r=n+\sum_i\sum_j m_{ij}p_i^{-j}$, where $n\in \ZZ$, the $p_i$'s run over the increasing 
sequence of all primes, the $m_{ij}$'s belong to $\{0,\dots,p_i-1\}$, and $j\in \NN^*$. 
Consequently, a generating subset of the additive group $\QQ$ is $\{1\} \cup\{1/p^j;\; p
\mbox{ a prime and } j\in \NN^*\}$. 
We let $z_1\in T_1$ be such that its class in $T_1/\langle1_{T_1} \rangle$ is $1$, 
and for every $n\in \ZZ$, we let $nz_1$ be the representative of the class 
$n$. If $x$ belongs to the class $1/p$ and $n\in \ZZ$, then 
there exists $k\in \ZZ$ such that 
$p(x+n1_{T_1})=px+pn1_{T_1}=z_1+k1_{T_1}+pn1_{T_1}$, 
so we see that we can choose $x$ in order to have 
$0\leq k\leq p-1$, we will denote this element by $k_p$. For $m \in \{ 2,\dots, p-1\}$, 
we let $mx$ be the representative of $m/p$. Assume that $j\geq 1$ and 
that the representative $x$ of the class $1/p^j$ is fixed, we let the representative 
of the class $1/p ^{j+1}$ be the element $y$ which satisfies  
$0\leq py-x\leq (p-1)1_{T_1}$, we let $k_{p^{j+1}}$ be the integer such that 
$k_{p^{j+1}}1_{T_1}=py-x$. 
For $m \in \{ 2,\dots, p-1\}$, we let $my$ be the representative of the class $m/p^{j+1}$. 
In the sum of any two rational numbers expressed as $r=n+\sum_i\sum_j m_{ij}p_i^{-j}$, 
the only interactions between the generators occur with $1/p^n$ and 
$1/p^{n-1}$ ($p$ a prime, $n\in \NN^*$). The mapping $\theta_1$ is the isomorphism between $T_1$ and $\QQ \overrightarrow{\times}_{\theta_1} \ZZ$. 
For every integers $m$ and $n$, we have $\theta_1(m,n)=0$, and $\theta_1$ is defined by 
induction in the same way as in 2) of Proposition \ref{propteta2}. We see that $\theta_1$ is 
uniquely determined by the $k_{p^j}$'s. We define $\theta_2$ in the same way. 
Now, the $k_{p^j}$'s are uniquely determined by the formulas 
$DD_{p^j,k}$, since $DD_{p^j,k}$ holds if and only if 
$k=k_p+k_{p^2}p+\cdots + k_{p^j} p^{j-1}$. It follows that $\theta_2=\theta_1$, and we have the isomorphism we were looking for. 
Note that if $p$ divides $z_1$, 
then the representative of the class of $1/p$ is the divisor of $z_1$, in this case, 
$k_p=0$, the same holds if $p^j$ divides $z_1$, hence the divisible hull of 
$\langle z_1\rangle$ in $T_1$ is contained in the set of representatives. 
\end{dem}
\subsection{Dense c-regular cyclically ordered abelian groups.}
\begin{num}{\bf Lemma.}\label{lm34} Let $p$ be a prime. If $G$ 
contains a $p$-torsion element, then $[p]\mbox{uw}(G)=[p]G$, otherwise, $[p]\mbox{uw}(G)=p[p]G$. 
\end{num}
\begin{dem}{\bf Proof.} 
 Let $x_1, \dots, x_n$ be elements of $\mbox{uw}(G)$ such that 
$\overline{x_1},\dots,\overline{x_n}$ are pairwise $p$-incongruent  
(within $G$), then $x_1,\dots,x_n$ are pairwise $p$-incongruent, 
hence $[p]\mbox{uw}(G)\geq [p]G$. If $[p]G$ is infinite, then the proposition is trivial. 
So we assume that $[p]G$ is finite, and let $\overline{x_1},\dots,\overline{x_n}$ be a 
maximal family of $p$-incongruent elements. 
In particular, for every $i\neq j$, $x_i-x_j \notin \langle z_G\rangle$. Let $x \in 
\mbox{uw}(G)$, there exists some $x_i$ such that $\overline{x}$ is $p$-congruent to 
$\overline{x_i}$, i.e. there exists $y\in \mbox{uw}(G)$ and $k\in \ZZ$ such that 
$x-x_i=py+kz_G$. If $z_G$ is divisible by $p$ (in other words, if $G$ contains 
a $p$-torsion element), this implies that $x$ and $x_i$ are 
$p$-congruent (within $\mbox{uw}(G)$), hence $x_1, \dots, x_n$ is a 
maximal $p$-incongruent family of elements of $\mbox{uw}(G)$. 
It follows: $[p]\mbox{uw}(G)= [p]G$. If $z_G$ is not $p$-divisible, 
then the elements $x_i+kz_G$, $1\leq i\leq n$, $0\leq k \leq p-1$, 
are pairwise $p$-incongruent. Let $x\in \mbox{uw}(G)$ and 
$x_i$ such that $\overline{x}=\overline{x_i}$, then, by what we did above, 
$x$ is $p$-congruent to one of the $x_i+kz_G$'s within $\mbox{uw}(G)$, which 
proves that the family we defined is maximal, and that $[p]\mbox{uw}(G)=p[p]G$. 
\end{dem}
\begin{num}{\bf Proposition.}\label{avprop36} 
Let $G_1$ and $G_2$ be two dense c-regular cyclically ordered groups 
such that $G_1$ is a subgroup of $G_2$. Then $G_1$ is an  
elementary substructure of $G_2$ if and only if $G_1$ is pure in $G_2$, and, 
for every prime $p$, $[p]G_1=[p]G_2$. 
\end{num} 
\begin{dem}{\bf Proof.} Trivially, if $G_1\prec G_2$, then $G_1$ is pure in $G_2$ and, 
for every prime $p$, $[p]G_1=[p]G_2$, because these are first order properties. 
Conversely, we know that  
$\mbox{uw}(G_1)$ embeds into $\mbox{uw}(G_2)$. First we prove that $\mbox{uw}(G_1)$ 
is pure in $\mbox{uw}(G_2)$. Let $(n,g)$ be an element of $G_1$ which is
divisible by some prime $p$ within $\mbox{uw}(G_2)$. Since $p(1,0)=(p,0)$, we can assume that  
$0\leq n<p$, hence a divisor of $(n,g)$ can be expressed as $(0,h)$, for some $h\in G_2$, 
and $ph=g$. Since $G_1$ is pure in $G_2$, it contains some element $h'$ such that $ph'=g$. 
If $h'\neq h$, then $h'h^{-1}$ is a $p$-torsion element within $G_2$, and 
since $G_1$ is pure in $G_2$, it contains a $p$-torsion element. 
We know that if a cyclically ordered group contains a $p$-torsion element, then  
its subgroup of $p$-torsion elements is cyclic of cardinal $p$ 
(to see this, look at the subgroups of $\KK$). Consequently $G_1$ contains all the 
$p$-torsion elements of $G_2$, and in particular, it contains $h'h^{-1}$, 
hence it contains $h$. Since any integer can be factored into prime powers, 
we deduce by induction that for every $m\in \NN^*$, 
if $m$ divides $(n,g)$ within $\mbox{uw}(G_2)$, then it divides 
$(n,g)$ within $\mbox{uw}(G_1)$ (note that the divisors are the same). 
We see also that for every prime $p$, $G_2$ contains a $p$-torsion element 
if and only if $G_1$ contains an $p$-torsion element, and by Lemma \ref{lm34}, 
we have $[p]\mbox{uw}(G_1)=[p]\mbox{uw}(G_2)$. By \cite{Za 61}, we have: 
$\mbox{uw}(G_1)\prec \mbox{uw}(G_2)$. Since $z_{G_2}=z_{G_1}\in \mbox{uw}(G_1)$, 
$\mbox{uw}(G_1)$ still remains an elementary substructure of $\mbox{uw}(G_2)$ in 
a language augmented with a predicate consisting of the constant $z_G$, and by Theorem 4.1 
of \cite{GLL 12}, we have: $G_1 \prec G_2$. 
\end{dem}
\begin{num}{\bf Proposition.}\label{prop317}  
Let $G_1$ and $G_2$ be two dense countable c-regular cyclically ordered groups 
having isomorphic torsion groups and  such that, for every prime $p$, 
$[p]G_1=[p]G_2$. Then $G_1$ and $G_2$ are elementary substructures of the same 
cyclically ordered group. 
\end{num}
\begin{dem}{\bf Proof.} Since $G_1$ and $G_2$ have the same prime invariants of Zakon  
and the same torsion groups, their linear parts have the same prime invariants of Zakon, 
hence, by 
Note 6.2 and Theorem 6.3 of \cite{Za 61}, they are elementary substructures of 
the same regular group $L$ where for every $u\in L$ there exists $n\in \NN^*$ such that 
$nu$ can be expressed as $g_1+g_2$, for some $g_1\in G_1$ and $g_2\in G_2$, 
the order being the lexicographic one. 
$\mbox{uw}(G_1)/l(G_1)$ and $\mbox{uw}(G_2)/l(G_2)$ embed into 
$\RR$, where the images of the classes of $z_{G_1}$ and $z_{G_2}$ are $1$. Let $R$ be the 
subgroup of $\RR$ generated by the images of 
$\mbox{uw}(G_1)/l(G_1)$ and $\mbox{uw}(G_2)/l(G_2)$, 
$R$ is divisible, because it is generated 
by two divisible subgroups. Let $\theta_1$ and $\theta_2$ be two mappings 
satisfying (*), (**) and (***) of Proposition \ref{propteta1}, and such that  
$\mbox{uw}(G_1)\simeq (\mbox{uw}(G_1)/l(G_1))\overrightarrow{\times}_{\theta_1}l(G_1)$ and 
$\mbox{uw}(G_2)\simeq (\mbox{uw}(G_2)/l(G_2))\overrightarrow{\times}_{\theta_2}l(G_2)$. 
By Proposition \ref{propteta2}, $\theta_1$ and $\theta_2$ extend to $R$. We set  
$\theta=(\theta_1,\theta_2)$, in the same way as in Proposition \ref{propteta1b}. 
According to the same proposition, we get a structure of linearly ordered group 
$R\overrightarrow{\times}_{\theta} L$ where $\mbox{uw}(G_1)$ 
and $\mbox{uw}(G_2)$ are cofinal subgroups. By properties of regular groups, 
since $L$ is regular and the quotient of $R\overrightarrow{\times}_{\theta} L$ 
by $L$ is divisible, the 
group $R\overrightarrow{\times}_{\theta}L$ 
is regular, and it has the same prime invariants of Zakon as $L$ 
and as $\mbox{uw}(G_1)$ and $\mbox{uw}(G_2)$. 
Since the theory of dense regular groups with fixed family of  prime invariants 
of Zakon is model complete, it follows that $\mbox{uw}(G_1)$ 
and $\mbox{uw}(G_2)$ are elementary substructures of 
$R\overrightarrow{\times}_{\theta}L$. Since 
these two subgroups contain $(1,0)$, they still remain elementary substructures 
in a language augmented with a 
predicate consiting of this element, and their wound-rounds $G_1$ and $G_2$ 
are elementary 
substructures of the wound-round associated to $R\overrightarrow{\times}_{\theta} L$. 
\end{dem}
\begin{num}{\bf Proposition.}\label{prop318} For any choice of the family  of prime invariants 
of Zakon and any subgroup of $\UU$, there exists a countable c-archimedean 
dense cyclically ordered group $G$ with family of 
prime invariants of Zakon and with torsion group so chosen. 
\end{num}
\begin{dem}{\bf Proof.} Let $T$ be the subroup of $\UU$, 
the unwound $\mbox{uw}(T)$ of $T$ is a subgroup of 
$\QQ$, and for every prime $p$, $n\in \NN$, the element $z_T$ is $p^n$-divisible 
if and only if 
$\zeta_{p^n}\in T$. In the same way as at the beginning of Section \ref{subsec31}, 
$[p]\mbox{uw}(T)=1$ if $\mbox{uw}(T)$ is $p$-divisible, and $[p]\mbox{uw}(T)=p$ otherwise. 
Let ${\cal Z}=\{p_i^{n_i} \mid i\in \NN^*\}$ be the given family of invariants of Zakon, 
where $(p_i)$ is the increasing sequence of all 
primes and the $n_i$'s belong to $\NN\cup \{\infty\}$. By Lemma \ref{lm33}, there  
exists a countable subgroup $M$ of $\RR$ such that for every $i\in \NN^*$, 
$[p_i]M=[p]\mbox{uw}(T)p_i^{n_i}$ and $\mbox{uw}(T)$ is a pure subgroup of $M$. 
We set $G=M/\langle z_T \rangle$. Since $\mbox{uw}(T)$ is pure in $M$, 
the groups $G$ and $T$ have the same torsion subgroup, which is  $T$, and by 
Lemma \ref{lm34}, we have, for every $i\in \NN^*$, $[p_i]G=p_i^{n_i}$. If $G$ is finite, hence 
it is contained in $\UU$, it suffices to consider the subgroup $H=\{e^{ir} \mid r\in \QQ\}$ of 
$\KK$, where $\QQ$ is the group of rational numbers, and for every 
prime $p$ we have 
$[p]H=1$. Hence $G\oplus H$ is infinite, so it is dense, it has the same torsion subgroup 
as $G$, and, for every prime $p$, $[p](G\oplus H)=[p]G \cdot [p]H=[p]G$. 
\end{dem}
\begin{num}{\bf Theorem.}\label{th313} 
There exists an infinite c-archimedean cyclically ordered group which is 
elementarily equivalent to $G$ if and only if $G$ is c-regular 
and dense. \\ 
Any two dense c-regular cyclically ordered groups are elementarily equivalent if and only if 
their torsion subgroups are isomorphic and they have the same family of prime invariants of 
Zakon. 
\end{num} 
\indent 
These conditions depend on the first order theory of $G$ by Lemma 
\ref{lm34} and because the number of $p$-torsion elements of $G$ depends on 
the first order theory of $G$. \\[2mm]  
\begin{dem}{\bf Proof.} If $G$ is elementarily equivalent to some c-archimedean infinite 
cyclically ordered group, then it shares with this group the first order properties, in particular 
being dense and c-regular. 
Now, assume that $G$ is c-regular and dense. By the theorem of 
L\"owenheim-Skolem, there exists a countable elementary substructure 
$G_1$ of $G$. By Proposition \ref{prop318}, there exists a countable c-archimedean dense 
group $G_2$ having the same torsion group and the same family 
of prime invariants of Zakon as $G_1$. By Proposition \ref{prop317}, there exists a 
cyclically ordered group $G'$ such that $G_1\prec G'$ and $G_2\prec G'$. 
Since $G_1\prec G$, we have $G_2\equiv G$. \\ 
\indent Let $G$ and $G'$ be c-regular dense cyclically ordered groups. If they are 
elementarily equivalent, then their torsion subgroups are isomorphic and they have 
the same family of prime invariants of Zakon. 
Assume now that the torsion subgroups of $G$ and $G'$ are isomorphic and that they 
have the same family of prime invariants of Zakon. 
Let $G_1\prec G$ and $G_1'\prec G'$ be countable. 
By Proposition \ref{prop317}, there exists a cyclically ordered group $G''$ such that 
$G_1\prec G''$ and $G_1'\prec G''$. It follows 
$G\equiv G'$. \end{dem}
\begin{num}{\bf Remark.} 
One can see, for example by Theorem \ref{prop4}, that every abelian c-divisible 
cyclically ordered group is c-regular, hence by Theorem \ref{th313} an abelian 
cyclically ordered group is c-divisible if and only if it is elementarily equivalent to $\UU$. 
\end{num}
\subsection{Discrete c-regular cyclically ordered abelian groups.}
\begin{num}{\bf Proposition.}\label{avprop38} 
Let $G_1$ and $G_2$ be two discrete c-regular cyclically ordered groups such that 
$G_1$ is a subgroup of $G_2$. Then $G_1$ is an 
elementary substructure of $G_2$ if and only if $G_1$ is pure in $G_2$ and 
the positive cones of $G_1$ and $G_2$ have the same first positive element.  
\end{num} 
\begin{dem}{\bf Proof.} Trivially, if $G_1\prec G_2$, then $G_1$ is pure in $G_2$, 
and their positive cones have the same first positive element, since it is definable. 
Conversely, in the same way as in the proof 
of Proposition \ref{avprop36}, we can show that $\mbox{uw}(G_1)$ is pure in 
$\mbox{uw}(G_2)$. Furthermore, $\mbox{uw}(G_1)$ and $\mbox{uw}(G_2)$ 
have the same lowest positive element, hence by \cite{Za 61}: 
$\mbox{uw}(G_1)\prec \mbox{uw}(G_2)$. Since $z_{G_2}=z_{G_1}\in \mbox{uw}(G_1)$, 
$\mbox{uw}(G_1)$ still remains an elementary substructure of $\mbox{uw}(G_2)$ in 
a language augmented with a predicate consisting of $z_G$, and by Theorem 4.1 
of \cite{GLL 12}, we have: $G_1 \prec G_2$. 
\end{dem}
\begin{num}{\bf Proposition.}\label{prop38} 
Assume that $G$ is c-regular discrete and is not 
c-archimedean. Then $G$ contains a discrete c-regular pure subgroup $H$ such that 
$K(H)=\UU$ and $l(H)\simeq \ZZ$, and which is an elementary substructure of $G$. 
\end{num}
\begin{dem}{\bf Proof.} Since $G$ is discrete, it contains a smallest nontrivial c-convex 
subgroup $\langle 1_G\rangle$, furthermore, $\langle 1_G\rangle$ is a pure subgroup 
of $\mbox{uw}(G)$, it is the smallest nontrivial convex 
subgroup of $\mbox{uw}(G)$. Since $G$ is c-regular, 
$\mbox{uw}(G)$ is regular and $\mbox{uw}(G)/\langle 1_G\rangle$ is divisible. 
Denote by $W$ the divisible hull of 
$\langle z_G\rangle\oplus \langle 1_G\rangle$ within $\mbox{uw}(G)$. 
$W$ is a pure subgroup of $\mbox{uw}(G)$ and 
$W/\langle 1_G\rangle$ is the divisible hull of $(\langle z_G\rangle\oplus \langle 
1_G\rangle)/\langle 1_G\rangle$ within $\mbox{uw}(G)/\langle 1_G\rangle$ 
(which is divisible), 
hence $W/\langle 1_G\rangle\simeq\QQ$. 
Furthermore, $W$ contains a unique proper convex subgroup: $\langle 1_G\rangle$, 
it follows that $W$ is a regular subgroup of $\mbox{uw}(G)$.
Since $W$ is pure in $\mbox{uw}(G)$ and it contains the same 
first positive element $1_G$, $W$ is an elementary substructure of 
$\mbox{uw}(G)$, and it still remains an elementary substructure in a language augmented 
with a 
predicate consisting of the element $z_G$. Set $H=W/\langle z_G\rangle$, by Theorem 4.1 of 
\cite{GLL 12}, $H$ is an elementary substructure of $G$. 
The linear part of $H$ is $\langle 1_G\rangle$, which is isomorphic to $\ZZ$, and 
$K(H)=(W/\langle z_G\rangle)/\langle 1_G\rangle\simeq \QQ/\langle 1\rangle =\UU$. 
\end{dem}
\indent The proof of this proposition shows that the first order theory of the 
infinite c-regular discrete cyclically ordered groups is equal to the first order theory 
of the abelian cyclically ordered groups whose linear part is isomorphic  
to $\ZZ$ and with quotient isomorphic to $\UU$, in other words, 
to the first order theory of the abelian linearly ordered groups $W$ which 
contain an unique convex subgroup, isomorphic to $\ZZ$, 
and such that the quotient group is isomorphic to the ordered group $\QQ$ and where 
the classe of $1$ is fixed. Note that $W$ is isomorphic to a subgroup of 
$\QQ\overrightarrow{\times} \QQ$ which contains  
$\ZZ\overrightarrow{\times}\{0\}$ and such that  
$\{0\}\overrightarrow{\times}\ZZ$ is its only convex subgroup. 
\begin{defi}{\bf Definition.} If $G$ is discrete and not c-archimedean, then 
the first positive element $1_G$ of $G$ is definable, 
we can assume that it lies in the language. For a prime $p$, $n\in \NN^*$ and $k\in \{0,\dots,p^n-1\}$, 
$D_{p^n,k}$ will be the formula: \\ 
$\exists x, \;R(0,x,2x,\dots, (p ^n-1)x)\wedge p^nx=k1_G$. 
\end{defi}
\begin{num}{\bf Remark.} 
Note that $G$ contains a $p^k$-torsion element 
whose class in $G/\langle1_G\rangle$ is a $p^k$-th root of $1$ 
if and only if $G$ satisfies the formula $D_{p^k,0}$, 
hence it contains an element of torsion $n=p_1^{k_1}\cdots 
p_j^{k_j}$, whose class is a $n$-th root of $1$ in $G/\langle 
1_G\rangle$ if and only if it satisfies the formulas 
$D_{p_1^{k_1},0},\dots ,D_{p_j^{k_j},0}$. In particular, $G$ 
contains a subgroup isomorphic to $\UU$ if and only if it 
satisfies the formulas $D_{p^k,0}$ for every prime number $p$ and $k\in \NN^*$. \end{num}
\begin{num}{\bf Lemma.}\label{lmtetadi1} Assume that $G$ 
is c-regular discrete and infinite. 
For every $p$ prime, $n\in \NN^*$ and $k\in \{0,\dots,p^n-1\}$, 
$G$ satisfies the formula $D_{p^n,k}$ if and only if $\mbox{uw}(G)$ satisfies  
the formula $DD_{p^n,k}$, where the fixed element of the language consists of $z_G$. 
\end{num} 
\begin{dem}{\bf Proof.} The formula $R(0,x,2x,\dots, (p ^n-1)x)\wedge p^nx=k1_G$ 
says that the $jx$'s where $1\leq j\leq p^n-1$ ``do not turn full circle'', 
but $p^nx$ ``turns full circle'', and is equal to $k1_G$. This is equivalent to saying that 
the image of $x$ in $\mbox{uw}(G)$ satisfies $p^nx=z_G+k1_G$. 
\end{dem}
\begin{num}{\bf Lemma.}\label{lmtetadi1b} 
Assume that $G$ is discrete, is not c-archimedean, and that 
$G/\langle 1_G\rangle$ is c-divisible. \\ 
1) For every $p$ prime and $n\in \NN^*$, there exists exactly one integer $k\in \{0,
\dots , p^n-1\}$ such that $D_{p^j,k}$ holds within $G$. \\ 
2) Let $p$ be a prime, $n\in \NN^*$, and $k\in \{0,\dots , p^n-1\}$ 
such that $D_{p^j,k}$ holds in $G$. We express $k$ as $a_0+a_1p+\cdots + 
a_{n-1}p^{n-1}$, where, for $1\leq j \leq n$, $a_j\in \{0,\dots, p-1\}$. Then for every 
$j\leq n$, and $k_j=a_0+\cdots+a_{j-1}p^{j-1}$, $D_{p^j,k_j}$ holds within $G$. \\ 
3) There exists an unique sequence $(\varphi_p)$, where $p$ runs over the increasing sequence 
of prime numbers and the $\varphi_p$'s are mappings from $\NN^*$ into $\{0,\dots,p-1\}$ 
such that for every prime $p$, every $n\in \NN^*$ and every $k\in \{0,\dots,p^n-1\}$, 
$D_{p^n,k}$ holds within $G$ if and only if $k=\varphi_p(1)+\varphi_p(2)p+
\cdots + \varphi_p(n)p^{n-1}$. 
\end{num} 
\begin{dem}{\bf Proof.} 1) and 2) follow from Lemma \ref{lmtetadi1}, Lemma 
\ref{lemteta5b}, and from the fact that $G$ is c-divisible if and only if 
$\mbox{uw}(G)$ is divisible. 3) is a consequence of 1) and 2). 
\end{dem}
\begin{defi}{\bf Definition.}\label{def322} Assume that $G$ is discrete, non-c-archimedean, and 
that $G/\langle 1_G\rangle$ is divisible. 
For every prime $p$, let $\varphi_p$ be the mapping from $\NN^*$ into $\{0,\dots,p-1\}$, 
defined in 3) of Lemma \ref{lmtetadi1b}, the sequence $(\varphi_p)$, will be called the  {\it 
cha\-rac\-te\-ris\-tic sequence} of $G$. 
\end{defi} 
\begin{num}{\bf Remark.}\label{rq322} 
Note that $G$ contains a subgroup which is isomorphic to $\UU$ 
if and only if for every prime $p$ the mapping $\varphi_p$ is the $0$ mapping. \end{num}
\begin{num}{\bf Proposition.}\label{protetadi6} 
Let $C_1$ and $C_2$ be two non-c-archimedean discrete cyclically ordered groups 
such that 
$C_1/\langle 1_{C_1}\rangle\simeq \UU\simeq C_2/\langle 1_{C_2}\rangle$. 
The following conditions are equivalent. \\ 
$C_1$ and $C_2$ are isomorphic \\ 
$C_1$ and $C_2$ satisfy the same formulas $D_{p^n,k}$ \\ 
$C_1 \equiv C_2$.  \\ 
The characteristic sequence of $C_1$ is equal to the characteristic sequence of $C_2$. 
\end{num} 
\begin{dem}{\bf Proof.} Trivially, if $C_1$ and $C_2$ are isomorphic, then they are 
elementarily equivalent, and if they are are elementarily  
equivalent  then they satisfy the same formulas $D_{p^n,k}$, 
since these are first order formulas. 
Now, assume that they satisfy the same formulas $D_{p^n,k}$. 
By Lemma \ref{lmtetadi1}, $\mbox{uw}(C_1)$ and $\mbox{uw}(C_2)$ 
satisfy the same formulas $DD_{p^n,k}$. Since 
$C_1/\langle 1_{C_1}\rangle\simeq \UU\simeq C_2/\langle 1_{C_2}\rangle$, we have  
$\mbox{uw}(C_1)/\langle 1_{C_1}\rangle\simeq \QQ\simeq \mbox{uw}(C_2)/\langle 1_{C_2}\rangle$, 
and by Proposition \ref{propteta6}, $\mbox{uw}(C_1)$ and $\mbox{uw}(C_2)$ 
are isomorphic in a language containing the element $z_C$. 
Consequently, $C_1$ and $C_2$ are isomorphic. Finally, by Lemma \ref{lmtetadi1b} 
the last condition is equivalent to the second one. 
\end{dem} 
\begin{num}{\bf Theorem.}\label{thm313}  
Any two non-c-archimedean c-regular discrete cyclically ordered groups are
elementarily equivalent if and only if they satisfy the same formulas $D_{p^n,k}$. 
\end{num} 
\begin{dem}{\bf Proof.} Since these formulas are first order formulas, 
if two c-regular discrete 
cyclically ordered groups are elementarily equivalent then they satisfy the same formulas $D_{p^n,k}$. Conversely, assume that $G$ and $G'$ are c-regular, discrete, and they 
satisfy the same formulas $D_{p^n,k}$. Let $H$ be the subgroup of $G$ and $H'$ 
be the subgroup of $G'$ 
defined in the proof  of Proposition \ref{prop38}. Since they are  elementary 
substructures, they satisfy the same formulas $D_{p^n,k}$, hence by 
P¡roposition \ref{protetadi6}, they are isomorphic. Hence $G$ and $G'$ have elementary  
substructures which are isomorphic, so they are elementarily equivalent. 
\end{dem}
\begin{num}{\bf Corollary.}\label{cor314}  
Any two c-regular discrete non-c-archimedean cyclically ordered groups 
are elementarily equivalent if and only if their characteristic sequences are equal. 
\end{num} 
\begin{dem}{\bf Proof.} This is a consequence of Theorem \ref{thm313} and of 
the equivalence between the second and the fourth condition of Proposition 
\ref{protetadi6}, since if $G$ is discrete, then $G$ is c-regular 
if and only if $G/\langle 1_G\rangle$ is divisible. 
\end{dem} 
\begin{num}{\bf Corollary.} If $G$ is c-regular and discrete, then 
the first order theory of $G$ is uniquely determined by 
the subgroup $H$ defined in Proposition \ref{prop38}. 
\end{num} 
\begin{dem}{\bf Proof.} By Theorem \ref{thm313}, 
the first order theory of $G$ is uniquely determined by 
the formulas $D_{p^n,k}$, now apply Proposition \ref{protetadi6}.  
\end{dem}
\begin{num}{\bf Proposition.}\label{prop426} For every sequence $(\varphi_p)$ of functions from 
$\NN^*$ into $\{0,\dots,p-1\}$, where $p$ runs over the increasing sequence of prime numbers, 
there exists a non-c-archimedean discrete c-regular cyclically ordered group $H$ whose 
characteristic sequence is $(\varphi_p)$. 
\end{num} 
\begin{dem}{\bf Proof.} Let $(p_n)_{n\in \Ni^*}$ be the increasing sequence of prime 
numbers. In the linearly ordered group $\QQ \overrightarrow{\times}\QQ$ 
we set $z_H=(1,0)$ and $1_H=(0,1)$. For $n\in \NN^*$ and $i\in \{1,\dots,n\}$, 
let $k_i=\varphi_{p_i}(1)+\varphi_{p_i}(2)p+\cdots+\varphi_{p_i}(n)p_i^{n-1}$. 
According to the chinese remainder theorem, there exists $k\in \{0,1,\dots,
(p_1\cdots p_n)^n-1\}$ such that for $i\in \{1,\dots,n\}$ $k$ is congruent to $k_i$ modulo 
$p_i^n$. Denote by $W_n$ the subgroup of $\QQ \overrightarrow{\times}\QQ$ generated 
by $\displaystyle{\left(\frac{1}{(p_1\cdots p_n)^n},\frac{k}{(p_1\cdots p_n)^n}\right)}$ 
and $1_H$. 
Hence 
$\displaystyle{(p_1\cdots p_{n-1})^n\left(\frac{1}{(p_1\cdots p_n)^n},\frac{k}{(p_1\cdots 
p_n)^n}\right)=\left(\frac{1}{p_n^n},\frac{k}{p_n^n}\right)}$ is congruent to 
$\displaystyle{\left(\frac{1}{p_n^n},\frac{k_n}{p_n^n}\right)}$ mo\-du\-lo $1_H$, 
the same holds for the other $i\in \{1,\dots, n\}$. It follows that $W_n$ satisfies the formulas 
$DD_{p_i^n,k_i}$, where $i\in \{1,\dots,n\}$. The sequence $(W_n)_{n\in \Ni^*}$ is an increasing sequence of subgroups of $\QQ \overrightarrow{\times}\QQ$, denote by $W$ the union of all the 
$W_n$'s. One can prove that $z_H$ is divisible modulo $1_H$ within $W$, hence 
$W/\langle z_H\rangle \simeq \QQ$, and that $W$ is discrete with first element $1_H$. 
Consequently, $H=W/\langle z_H\rangle$ is a discrete c-regular subgroup such that 
$K(H)\simeq \UU$, $l(H)\simeq \ZZ$, the sequence $(\varphi_p)$ is its characteristic 
sequence. 
\end{dem} 
\indent In Remark \ref{rq322} we characterized  
the discrete c-regular cyclically ordered groups whose torsion 
subgroups are isomorphic to $\UU$, we can also characterize those which are divisible. 
The following result proves that this case can be seen as the opposite case. 
\begin{num}{\bf Proposition.}\label{prop327d}  
Let $G$ be a discrete c-regular cyclically ordered group, 
$p$ be a prime and $\varphi_p$ be the mapping from $\NN^*$ to $\{0,\dots,p-1\}$ defined 
in \ref{def322}. The following conditions are equivalent. \\ 
1) $G$ is $p$-divisible. \\ 
2) $1_G$ is $p$-divisible. \\ 
3) $\varphi_p(1)\neq 0$. \\
4) $G$ doesn't contain any $p$-torsion element. 
\end{num} 
\begin{dem}{\bf Proof.} We have trivially: 1) $\Rightarrow$ 2) and 4) $\Rightarrow$ 3). \\ 
\indent Assume that $\varphi_p(1)\neq 0$. Hence, for every 
$n\in \NN^*$, $p$ doesn't divide $k_{p^n}$. By 
Bezout identity, there exist two integers $u$ and $v$ such that $up^n+vk_{p^n}=1$. 
We know that there  exists $x$ in $G$ such that $p^nx=k_{p^n}1_G$, 
hence $p^n(vx)=vk_{p^n}1_G=
(1-up ^n)1_G$, and $p^n(vx+u)1_G=1_G$, 
wich proves that $1_G$ is $p$-divisible. Let $y\in G$, 
since $G$ is c-regular, the class of $y$ modulo $\langle 1_G\rangle$ is 
$p^n$-divisible, consequently there exists $z\in G$ such that $p^nz-y\in \langle 1_G\rangle$, hence there exists an integer $k$ such that $p^nz-k1_G=y$. Since $1_G$ is $p^n$-divisible, 
it follows that 
$y$ is $p^n$-divisible. We proved: 3) $\Rightarrow$ 2) $\Rightarrow$ 1). \\ 
\indent Assume that $G$ contains a $p$-torsion element $y$, then $2y,\dots, 
(p-1)y$ are $p$-torsion elements. By taking the lowest element of the set 
$\{y, 2y,\dots,(p-1)y\}$ ordered by $<_0$ instead of $y$, we can assume:  
$R(0,y,2y,\dots,(p-1)y)$, 
and $py=0$, so, by the definition of $k_p$, we have: $k_p=0$, that is, $\varphi_p(1)=0$. 
We proved 3) $\Rightarrow$ 4). 
Assume furthermore that $1_G$ is $p$-divisible, hence there exists $z$ such that $pz=1_G$. 
Let $j\in \{0,\dots, p-1\}$ be such that $R(jy,z,(j+1)y)$ holds ($z$ is not a multiple of $y$ 
because it is not $p$-torsion), by setting $z'=z-jy$ we have $R(0,z',y)$, wich implies 
$R(0,z',2z',\dots, pz')$, hence  
$R(0,z',1_G)$, a contradiction, because $1_G$ is minimal. It follows that 
$1_G$ is not $p$-divisible, hence $G$ is not $p$-divisible. Consequently, 2) $\Rightarrow$ 4), 
and the proposition is proved. 
\indent 
\end{dem}
\begin{num}{\bf Proposition.}
Let $G$ be a non-c-archimedean cyclically ordered group, then  $G$ 
is c-regular discrete divisible if and only if there exists a discrete cyclic order $R$ on the group 
$\QQ$ such that $(\QQ,R)$ is an elementary subextension of $G$. 
\end{num}
\begin{dem}{\bf Proof.} First assume that $G$ is c-regular divisible and discrete, 
then according to Proposition \ref{prop38}, there exists an elementary subextension $H$
of $G$ such that $l(H)\simeq \ZZ$ and $K(H)\simeq \UU$. $\QQ \cdot 1_H$ is a subgroup of $H$, 
because $H$ is divisible. Let $x\in H$, since $H/\langle 1_H\rangle\simeq\UU$, there exist 
integers $n$ and $k$ such that $nx=k1_H$. Since $H$ is divisible and torsion-free, 
(by Proposition \ref{prop327d}), we have: $x=\frac{k}{n} 1_H$, which proves: 
$H=\QQ \cdot 1_H$. In order to prove the other implication, it suffices to show that if 
$H=(\QQ,R)$ where $R$ is a discrete cyclic order, then $H$ is c-regular. 
Now, one can prove that $l(H)=\langle 1_H\rangle$, hence it is regular, and that 
$K(H)\simeq \UU$, consequently, 
by Theorem \ref{prop2}, $H$ is c-regular. 
\end{dem} 
\begin{defi}{\bf Definitions.} For every prime $p$, let $\varphi_p$ 
be a mapping from 
$\NN^*$ into $\{0,\dots,p-1\}$. For every prime $p$ and every $n\in \NN^*$, 
denote by 
$N_{p^n,\varphi_p}$ the set $p^n\NN^*-(\sum_{k=1}^n p^{k-1}\varphi_p(k)$). 
The set of the $N_{p^n,\varphi_p}$'s will be called the 
{\it family of subsets of} 
$\NN^*$ {\it characteristic of} $(\varphi_p)$. The family of subsets of  
$\NN^*$ characteristic of the characteristic sequence of $G$ will be called the  
{\it family of subsets of} $\NN^*$ {\it characteristic of} $G$. 
\end{defi} 
\begin{num}{\bf Proposition.}\label{prop315} \\ 
1) For every ultrafilter $U$ on $\NN^*$, there exists one and only one sequence 
$(\varphi_p)$, where $p$ runs over the increasing sequence of all primes
and $\varphi_p$ is a mapping from $\NN^*$ into $\{0,\dots,p-1\}$, 
such that $U$ contains the family of subsets of $\NN^*$ characteristic of 
$(\varphi_p)$. \\ 
2) For every $(\varphi_p)$, where $p$ runs over the increasing sequence of prime numbers 
and $\varphi_p$ is a mapping from $\NN^*$ into $\{0,\dots,p-1\}$, 
there exists a non principal ultrafilter $U$ on $\NN^*$ containing the family of subsets of 
$\NN^*$ characteristic of $(\varphi_p)$. 
\end{num}
\begin{dem}{\bf Proof.} 1) By a property of ultrafilters, if $A\in U$ and 
$A=A_1\cup \cdots \cup A_n$ is a finite partition of $A$, then $U$ contains 
exactly one of the subsets $A_k$, $1\leq k\leq n$. Now, for $p$ prime,  
$n\in \NN^*$, and $a_1,\dots, a_n$ in $\{ 0,\dots, p-1\}$, 
$p^{n+1}\NN^*-(\sum_{k=1}^na_kp^{k-1}+jp^n)$, ($0\leq j\leq p-1$) is a finite  
partition of $p^n\NN^*-(\sum_{k=1}^na_kp^{k- 1})$, the result follows by induction. \\ 
\indent 2) Let $p$ be fixed and $0<n_1<n_2<\cdots < n_k$, then $N_{p^{n_1},\varphi_p} 
\cap \cdots \cap N_{p^{n_k},\varphi_p}=N_{p^{n_k},\varphi_p}$, and if $1<p_1<p_2<
\cdots<p_k$ are prime, then by the chinese remainder theorem  
$N_{p_1^{n_1},\varphi_{p_1}} \cap \cdots \cap N_{p_k^{n_1},\varphi_{p_k}}$ is infinite. 
It follows that the intersection of every finite family of sets $N_{p^n,\varphi_p}$ has 
infinite cardinal. We join the elements of the filter of cofinite subsets to this family, then 
the property of non-empty intersection still remains satisfied, hence there exists an 
ultrafilter containing all the sets $N_{p^n,\varphi_p}$ and the filter of cofinite subsets. 
This ultrafilter is not principal, because it contains all the cofinite subsets. 
\end{dem}
\begin{defi}{\bf Definition.} Let $U$ be an ultrafilter on $\NN^*$, the  
family of subsets defined in 1) of Proposition \ref{prop315} will be called  
{\it the family of subsets of} 
$\NN^*$ {\it defined by} $U$. 
\end{defi} 
\begin{num}{\bf Theorem.}\label{thm316} \\ 
1) Let $U$ be an non principal ultrafilter on $\NN^*$, $C$ be the ultraproduct of the 
cyclically ordered groups $\ZZ/n\ZZ$ modulo $U$, $p$ be a prime, $n\in \NN^*$ and 
$k\in \{0,\dots,p^n-1\}$. Then $C$ satisfies the formula $D_{p^n,k}$ if and only if 
$p^n\NN^* -k\in U$. \\ 
2) Let $U$ be a non principal ultrafilter on $\NN^*$, $C$ be the ultraproduct of the 
cyclically ordered groups $\ZZ/n\ZZ$ modulo $U$. The family of subsets of $\NN^*$ 
defined by $U$ is equal to the family of subsets of $\NN^*$ characteristic 
of $C$. \\ 
3) Let $U_1$ and $U_2$ be two non principal ultrafilters on $\NN^*$, $C_1$ (resp. $C_2$) 
be the ultraproduct of the cyclically ordered groups $\ZZ/n\ZZ$ modulo $U_1$ (resp. $U_2$). 
Then $C_1 \equiv C_2$ if and only if the family of subsets of $\NN^*$ 
defined by $U_1$ is the same as the family of subsets of $\NN^*$ 
defined by $U_2$. 
\end{num}
\begin{dem}{\bf Proof.} 1) $C$ satisfies $D_{p^n,k}$ if and only if 
the set of all integers $m$ such that $\ZZ/m\ZZ$ satisfies $D_{p ^n,k}$ belongs to $U$. 
Let $m\in \NN^*$, there exists some $j\in \{0,\dots, p^n-1\}$ such that $m \in p^n\NN^*-j$, 
then there exists $x\in \NN^*$ such that $m=xp^n-j$, and we have 
$p^nx=m+j$, and $0<x<2x<\cdots <(p^n-1)x<m$. In $\ZZ/m\ZZ$, this is equivalent to $R(0,x,2x,\dots, (p^n-1)x)$ and $p^nx=j$, hence $\ZZ/m\ZZ$ satisfies 
$D_{p^n,j}$, and for $j'\neq j$, $0\leq j'\leq p^n-1$, $\ZZ/m\ZZ$ does not satisfy $D_{p^n,j'}$, 
since every discrete c-regular cyclically ordered group satisfies exactly 
one relation $D_{p^n,j}$, for fixed $p$ and $n$. It follows that $\ZZ/m\ZZ$ satisfies 
$D_{p^n,k}$ if and only if $m\in p^n\NN^*-k$. Consequently, 
$C$ satisfies the formula $D_{p^n,k}$ if and only if 
$p^n\NN^* -k\in U$. \\ 
\indent 2) Follows from 1) and from Lemma \ref{lmtetadi1b}. \\
\indent 3) Follows from 2) and from Corollary \ref{cor314}. 
\end{dem} 
\begin{num}{\bf Corollary.} 1) If $G$ is infinite, c-regular and discrete, and $U$ is a 
non principal ultrafilter on $\NN^*$, then $G$ is elementarily equivalent 
to the ultraproduct of the cyclically ordered groups $\ZZ/n\ZZ$ modulo $U$ if and 
only if the family of subsets of $\NN^*$ characteristic of $G$ is 
the same as the family of subsets of $\NN^*$ defined by $U$. \\ 
2) $G$ is c-regular and discrete if and only if there exists an ultrafilter on $\NN^*$ such 
that $G$ is elementarily equivalent to the ultraproduct of the 
cyclically ordered groups $\ZZ/n\ZZ$ modulo $U$. 
\end{num}
\begin{dem}{\bf Proof.} 1) Follows from Theorem \ref{thm316} and from Corollary 
\ref{cor314}. \\ 
2) If there exists an ultrafilter on $\NN^*$ such that $G$ is elementarily equivalent 
to the ultraproduct of the cyclically ordered groups $\ZZ/n\ZZ$ modulo $U$, then 
it is c-regular and discrete, because these are first order properties. 
Conversely, assume that $G$ is c-regular and discrete. If $G$ is finite, then 
there exists $n_0\in \NN^*$ such that $G=\ZZ/n_0\ZZ$, then $G$ is isomorphic to 
 the ultraproduct of the cyclically ordered groups $\ZZ/n\ZZ$ modulo the principal 
ultrafilter  generated by $\{n_0\}$. If $G$ is infinite, the result follows from 1). 
\end{dem} 
\begin{num}{\bf Proposition.} For every subgroup $S$ of $\UU$, there exists an ultraproduct 
of finite cyclic groups whose torsion subgroup is $S$. 
\end{num}
\begin{dem}{\bf Proof.} This proposition is a consequence of Theorem 
\ref{thm316}. 
\end{dem}
\begin{num}{\bf Theorem.} 1) The class of all discrete c-regular cyclically ordered abelian groups 
is the smallest elementary class which contains all finite cyclic groups. \\ 
2) The class of all discrete regular cyclically ordered abelian groups 
is the smallest elementary class which contains all cyclic groups. 
If $G$ belongs to this class, then either $G$ is a linearly cyclically ordered group 
which is elementarily equivalent to $\ZZ$, or $G$ is c-regular and elementarily equivalent to 
an ultraproduct of finite cyclic groups. 
\end{num}
\begin{dem}{\bf Proof.} Follows from what we just proved, from Theorem \ref{prop2b} 
and from the properties of regular linearly ordered abelian groups.  
\end{dem} 
\begin{num}{\bf Proposition.}\label{propteta4}
Every linearly ordered abelian group $T$, having a smallest proper convex subgroup, 
embeds into a minimal regular group $T'$ which is contained into its divisible hull. If 
$T$ is dense, then $T'$ is dense, if $T$ is discrete, then $T'$ is discrete. In 
particular, every discrete linearly ordered abelian group embeds into a 
minimal discrete regular group which is contained in its divisible hull. 
\end{num}
\begin{dem}{\bf Proof.} Denote by $C$ the smallest proper convex subgroup of $T$, 
then there exists a cocycle $\theta$ such that $T\simeq (T/C)
\overrightarrow{\times}_{\theta} C$. We embed $(T/C)$ into its divisible hull and 
we extend $T$ in the same way as in Proposition \ref{propteta2} 
(note that this embedding need not be unique) 
the details are left to the reader. \end{dem}
\begin{num}{\bf Proposition.} Every dense cyclically ordered group having a 
smallest proper c-convex subgroup embeds into a minimal c-regular 
cyclically ordered group which is contained in its divisible hull. 
\end{num}
\begin{dem}{\bf Proof.} By Proposition \ref{propteta4}, $\mbox{uw}(G)$ 
embeds into a minimal regular group contained in its divisible hull, we 
consider the wound-round associated to its divisible hull and $z_G$. 
By Theorem \ref{prop4}, this 
extension is minimal. 
\end{dem}
\begin{num}{\bf Proposition.}\label{propos316} 
Assume that $G$ is $\omega_1$-saturated and $K(G)$ 
is infinite. \\ 
1) Every class of $G$ modulo $l(G)$ contains a divisible element. \\ 
2) Every class of $\mbox{uw}(G)$ modulo $l(G)$ contains a divisible element. \\ 
3) Let $p$ be a prime such that $G$ is not $p$-divisible, then every class of $G$ 
modulo $l(G)$ contains an element which is not $p$-divisible and which is $n$-divisible 
for all integers $n$ such that $n$ and $p$ are coprime. \\
4)  Let $p$ be a prime such that $G$ is not $p$-divisible, then every class of 
$\mbox{uw}(G)$ modulo $l(G)$ contains an element which is not $p$-divisible and 
is $n$-divisible for all integers $n$ such that $n$ and $p$ are coprime. 
\end{num}
\begin{dem}{\bf Proof.} 1)  Note that by Propositon 6.3 a) of \cite{GLL 12}, $K(G)\simeq\KK$, 
since $G$ is $\omega_1$-saturated and $K(G)$ is infinite. 
Let $x_0\in G$, assume first that $\overline{x_0}\notin \UU$, say 
$\overline{x_0}=e^{i\alpha}$ for some irrational element $\alpha \in [0,1[$. 
For every $n\in \NN^*$, let $m_n$ be the integer such that 
$2\pi\frac{m_n}{n} \leq \alpha< 2\pi\frac{m_n+1}{n}$. 
For any $n\in \NN^*$, $x\in G$ and $\theta \in [0,1[$ such that $\bar{x}=e^{i\theta}$, 
we have: $\frac{2\pi}{n} < \theta< \frac{2\pi}{n-1} \Rightarrow R(0,x,2x,\dots,(n-1)x)\& \neg R(0,x,2x,\dots,(n-1)x,nx))\Rightarrow \frac{2\pi}{n} \leq \theta\leq \frac{2\pi}{n-1}$, 
hence one can see that 
$2\pi\frac{m_n}{n}< \theta< 2\pi\frac{m_n+1}{n}$ implies 
$$[\forall y,(R(0,y,2y,\dots,(n-1)y)\& \neg R(0,y,2y,\dots,(n-1)y,ny)) \Rightarrow 
R(0,m_ny,x)]$$ 
$$[\exists y,R(0,y,2y,\dots,(n-1)y)\&\neg R(0,y,2y,\dots,(n-1)y,ny) 
\&R(0,x,(m_n+1)y)]$$ 
which implies $2\pi\frac{m_n}{n} \leq \theta\leq 2\pi\frac{m_n+1}{n}$. 
In any case, $\bar{x}=\bar{x_0}$ if and only if $x$ satisfies all of these first order formulas. 
The element $x$ is divisible if and only if 
$\exists y, \; ny=x$ holds for every $n\in \NN^*$. We get a countable type. 
Take a finite subset of this type, and let $m$ be lcm of all the $n$'s that appear in 
some formula $\exists y,\; ny=x$. Then pick an element  
$y$ of $\frac{1}{m}\overline{x_0}$, and $x=my$, $x$ satisfies this finite set of formulas. 
By $\omega_1$-saturation, the type is satisfied, and $\overline{x_0}$ contains 
a divisible element. \\ 
\indent Now assume that $\overline{x_0}\in \UU$, say $\overline{x_0}=
e^{2ik\pi/n}$, where $k$ and $n$ are coprime, $k<n$. In $\KK$, 
$\overline{x_0}$ is characterized by $\overline{x_0}\neq 0, \; 
2\overline{x_0}\neq 0,\dots, (n-1)\overline{x_0}\neq 0, n\overline{x_0}=0$, and 
for $2\leq l\leq n-1$, one of the relations $R(0,\overline{x_0},l\overline{x_0})$ 
or $R(0,l\overline{x_0},\overline{x_0})$. 
We are going to find sufficient conditions in order that these formulas be satisfied. 
For $n\overline{x}=0$, i.e. 
$x\in l(G)$, we take all the formulas $R(0,nx,knx)$, for $k\geq 2$. For 
$l\overline{x}\neq 0$, i.e. $lx\notin l(G)$, we take the smallest integer $k_l$ such that 
$R(0,k_l l\overline{x_0},l\overline{x_0})$ holds, and we take the formula 
$R(0,k_l lx,lx)$. $x$ is divisible if and only if $\exists y,\; qy=x$ holds for every 
$q\in \NN^*$. This type is countable. 
Take a finite subset of this type, let $m$ be the lcm of the $q$'s such that 
$\exists y,\; qy=x$ appear, and $y\in \frac{1}{m}\overline{x_0}$. Then 
$mny\in l(G)$. If $mny$ does not belong to the positive cone of $G$, we take 
$y'=y-mny$, then $mny'=mny-mnmny$ belongs to positive cone of $G$. Set   
$x=my'$, and the finite subset of formulas is satisfied. By $\omega_1$-saturation, 
the type is satisfied. \\ 
2) The elements of $\mbox{uw}(G)$ can be expressed as $(n,x)$, for some $n\in \ZZ$ and 
$x\in G$. In order to prove that $(n,x)$ is divisible in $\mbox{uw}(G)$, 
we can assume that 
$n\geq 0$, and it is sufficient to prove that $(n,x)$ is divisible by every integer $m\geq n$. 
In this case, the divisor is some $(0,y)$ with $y\in G$. Then $m(0,y)=(n,x)$ is 
equivalent to: 
the sequence $y,2y,\dots, my$ ``travels $n$ times around'' and furthermore 
$my=x$. $(k+1)y$ ``being in the round which follows the round'' of $y$ is equivalent to: 
$R(0,(k+1)y,ky)$. Since $m$ and the class of $x$ are fixed, we know where are the jumps, 
hence we have a succession of $R(0,(k+1)y,ky)$ and $R(0,ky,(k+1)y)$ 
which is well-defined. The class of $(n,x)$ modulo $l(G)$ containing a divisible element 
can be rendered by: there exists $z$ such that $z-x\in l(G)$ and $(n,z)$ is divisible. 
Finally $z-x\in l(G)$ is equivalent to an infinite set of formulas $R(0,z-x,q(z-x))$, 
for $q\geq 2$. This type is countable. 
In order to get an element which satisfies a finite subset of this type, we consider the lcm 
of all the $m$'s which appear in this subset (we still denote it by $m$), we pick  
$(0,t)$ in the class of $\frac{1}{m}(n,x)$ modulo $l(G)$ (which exists since 
$\mbox{uw}(G)/l(G)$ is divisible), then $z=mt$, satisfies the finite subset of formulas. 
The only formulas which appear belong to the language of cyclically ordered groups. 
By $\omega_1$-saturation, the type is satisfied, which gives us the required divisible 
element. \\
3)  and 4) Let  $x_0$ be an element of $G$ which is not $p$-divisible, and we consider the 
countable type $\forall y,\; py\neq x$, and $\exists y,\; ny=x$, 
for every positive integer $n$ prime to $p$. Consider a finite subset of this type, 
and let $m$ be the lcm of all the $n$'s prime to $p$ which appear. Then $x=mx_0$ 
satisfies this finite subset. By $\omega_1$-saturation, $G$ contains an element 
$x$ which is not divisible by $p$ and which is divisible by every integer which is prime to $p$. 
By 1), there exists $y$ which is divisible and such that $\bar{y}=\bar{x}$, hence 
$x-y\in l(G)$ is not $p$-divisible and is divisible by every integer which is prime to $p$. 
Now, we add to $x-y$ some divisible element of a fixed class of  
$G$ or of $\mbox{uw}(G)$. 
\end{dem}
\begin{num}{\bf Proposition.}\label{prop36} Assume that $G$ is c-regular, dense and 
$\omega_1$-saturated. Then $G$ contains a countable elementary substructure 
which is a c-ar\-chi\-me\-dean group.  
\end{num} 
\begin{dem}{\bf Proof.} 
Denote by $H$ the divisible hull of $\langle z_G\rangle$ in 
$\mbox{uw}(G)$. Since $H$ is a pure subgroup of $\mbox{uw}(G)$, 
we have for every prime $p$: 
$[p]H\leq [p]\mbox{uw}(G)$. By Lemma \ref{lm33}, there exists a countable subgroup 
$M$ of $\RR$ containing a pure subgroup isomorphic to $H$ and such that for 
every prime $p$ we have $[p]M=[p]\mbox{uw}(G)$. 
We assume that $\mbox{uw}(G)/l(G)=\RR$
by letting $1$ being the class of $z_G$; $M$ is a subgroup of $\mbox{uw}(G)/l(G)$ 
that we embed into $\mbox{uw}(G)$ in the following way. 
The image of $1$ is $z_G$. Every $e^{k_{n_j}}$ is 
a class modulo $l(G)$, its image is some representative which is not 
divisible by $p_n$ and is divisible by every integer which is prime to $p_n$. 
This element exists by Proposition \ref{propos316}. So, 
every $e^{k_{n_j}}$ satisfies the same divisibility relations within $M$ as its image 
within $\mbox{uw}(G)$. By the definition of $M$ its image is a pure sbgroup of 
$\mbox{uw}(G)$. Now, this image is indeed an ordered subgroup, 
because the images of the $e^{k_{n_j}}$'s are the corresponding classes. 
$M$ is a pure subgroup of $\mbox{uw}(G)$ such that, for every prime $p$,  
$[p]M=[p]\mbox{uw}(G)$, furthermore $M$ is regular since it is archimedean. It follows that 
$M\prec \mbox{uw}(G)$. Consequently, $M/\langle z_G\rangle$ is a countable 
elementary substructure of $G$, which is a c-archimedean subgroup. 
\end{dem}
\begin{num}{\bf Proposition.} 
Every discrete cyclically ordered group embeds into a c-regular discrete 
cyclically ordered group contained in its divisible hull. 
\end{num}
\begin{dem}{\bf Proof.}  By Proposition \ref{propteta4}, $\mbox{uw}(G)$ 
embeds into a minimal regular group contained in its divisible hull, we 
consider the wound-round associated to this divisible hull and $z_G$. 
The minimality follows from Theorem \ref{prop4}. 
\end{dem}
\section{With or without the cyclic order predicate?}
\indent 
It seems to be natural and usefull for studying $\RR$ to consider its linear order, what about 
the study of $\KK$ and its cyclic order? More generally, to determine the contribution 
of the cyclic order in the study of the c-regular cyclically ordered groups, we focus on the 
theory of those groups either in the language for groups without a predicate for the cyclic 
order or with the cyclic order predicate. In the strict language of groups,  
we know that any two abelian groups 
are elementarily equivalent if and only if they have the same invariants of Szmielew 
(cf. \cite{Sz 55}). Here the torsion subgroups of our groups embed into $\UU$, hence the 
Szmielew invariants reduce to the Zakon invariants and the torsion subgroup. 
We can reformulate the second part of Theorem \ref{th313} as follows. 
\begin{num}{\bf Proposition.} Any two dense c-regular cyclically ordered abelian groups are 
elementarily equivalent 
if and only if they are elementarily equivalent groups. \end{num}
\indent 
The same holds for any two dense regular linearly ordered groups. \\ 
\indent With respect to the discrete c-regular cyclically ordered groups, 
all groups having a fixed torsion subgroup have the same Szmielew 
invariants, it follows: 
\begin{num}{\bf Proposition.} Any two infinite discrete c-regular cyclically ordered abelian 
groups having a fixed torsion subgroup are elementarily equivalent in the language of groups. \\ 
Two infinite discrete c-regular cyclically ordered abelian 
groups having a fixed torsion subgroup need not be elementarily equivalent in the language of 
cyclically ordered groups.
\end{num}
Indeed, by Theorem \ref{thm313}, 
they are elementarily equivalent in the language of cyclically ordered groups if and 
only if they satisfy the same formulas $D_{n,k}$. Now, by Lemma \ref{lmtetadi1b} and 
Proposition \ref{prop426}, they need not satisfy the same formulas $D_{n,k}$. \\[2mm]
\indent In the case of linearly ordered groups, we have the following. 
\begin{num}{\bf Proposition.}\label{prop53} If $G$ is a torsion-free abelian group, then 
there exists a structure of regular linearly ordered group on $G$ with a smallest 
convex proper subgroup if and only if for every prime $p$ the cardinal of every
maximal family of pairwise $p$-incongruent elements is 
at most equal to the cardinal of $\RR$. Moreover, 
there exists a structure of discrete regular linearly ordered group on $G$ 
if and only if, for every $p$ prime, $[p]G=1$ and there exists an  
element which is not divisible by any prime $p$. \end{num}
\begin{dem}{\bf Proof.} Assume that for every prime $p$ the cardinal of every
maximal family of pairwise $p$-incongruent elements is 
at most equal to the cardinal of $\RR$ and for each prime $p$ let $I_p$ be 
such a family. Denote by $H$ the divisible hull   
within $G$, of the subgroup generated by the union of all the subsets $I_p$, 
where $p$ run over the set of all prime numbers. 
The cardinal of $H$ is at most equal to the cardinal of $\RR$, hence the cardinal 
of a basis of the $\QQ$-vector space $\tilde{H}$, where $\tilde{H}$ is the divisible hull of 
$H$, is at most equal to the cardinal 
of $\RR$. We embed this basis into a rationally independent subset of 
$\RR$, so we get an embedding of the group $\tilde{H}$ into $\RR$, hence  
we have an embedding of $H$ into $\RR$. This embedding gives rise to a linear order on $H$, 
such that $H$ is an archimedean group. Now, the 
quotient group $G/H$ is torsion-free, hence there is a compatible linear order on $G/H$. 
Let $H'$ be a set of representatives of the equivalence classes of $G$ modulo $H$,   
we deduce a one-to-one mapping from $G$ onto $H'\times H$. 
The lexicographic order on $H'\times H$ induces a linear order on $G$, and in the same way as 
in \cite{Ja 54}, one can check that this order is compatible. Furthermore, 
$G/H$ is divisible and $H$ is archimedean, hence $G$ is regular. Conversely, 
let $C$ the smallest convex proper subgroup of $G$. Since $G/C$ is 
divisible, for every prime $p$ there is a maximal family $I_p$ of pairwise $p$-incongruent 
elements which is contained in $C$. Now, $C$ is archimedean, hence 
$\mbox{card}(C)\leq\mbox{card}({\RR})$, so $\mbox{card}(I_p)\leq\mbox{card}({\RR})$. 
By \cite[Corollary 1.8]{Za 61} every maximal family of pairwise $p$-incongruent elements 
has the same cardinal, hence the cardinal of every
maximal family of pairwise $p$-incongruent elements is 
at most equal to the cardinal of $\RR$. \\ 
\indent If for every prime $p$ we have $[p]G=1$ and there exists an element $g_0$ 
which is not divisible by any prime, we let $C$ be the subgroup generated by $g_0$, it is 
a pure subgroup and $G/C$ is divisible. We conclude in the same way as above. 
The converse is trivial. 
\end{dem} 
\indent Now assume that $G$ is a torsion-free abelian group such that 
there exists a prime $p$ and a maximal family $I_p$ of elements which 
are pairwise $p$-incongruent and of cardinal greater than the cardinal of 
$\RR$ and, for every prime $q\neq p$, $G$ is $q$-divisible. Let $H$ be the 
divisible hull within $G$ of the subgroup generated by $I_p$. 
Let $C$ be a pure subgroup of $H$ distinct from $H$. Hence one of the 
elements of the family $I_p$ does not belong to $C$, hence $H/C$ is not divisible. 
Since the cardinal of $H$ is greater than the cardinal of $\RR$, there cannot 
exist a structure of archimedean linearly ordered group on $H$. 
Hence for every linear order $H$ contains a proper convex subgroup $C$, now 
$H/C$ is not divisible, 
so, there does not exist any structure of regular linearly ordered group on $H$. 
For the same reason, 
there does not exist any structure of c-regular cyclically ordered group. \\[2mm] 
\indent Turning to the cyclically ordered case we get: 
\begin{num}{\bf Proposition.}\label{prop54} If $G$ is a torsion-free abelian group, then 
there exists a structure of c-regular cyclically ordered group on $G$ with a smallest 
convex proper subgroup if and only if for every prime $p$ the cardinal of every
maximal family of pairwise $p$-incongruent elements is 
at most equal to the cardinal of $\RR$. Moreover, 
there exists a structure of discrete c-regular cyclically ordered group on $G$ 
if and only if, for every $p$ prime, $[p]G=1$ and there exists an  
element which is not divisible by any prime $p$. \end{num}
\begin{dem}{\bf Proof.}
Leg $G$ be an abelian group such that 
for every prime $p$ the cardinal of every
maximal family of pairwise $p$-incongruent elements is 
at most equal to the cardinal of $\RR$ and for each prime $p$ let $I_p$ be 
such a family. 
Let $H$ be the divisible hull within $G$ of the subgroup generated by the union of all 
the $I_p$'s (hence the cardinal of $H$ is at most equal to the cardinal of $\RR$), 
let $A$ be a $\QQ$-vectorial subspace of the divisible group $G/H$, 
of cardinal at most equal to the cardinal 
of $\RR$, and let $B$ be a subspace such that 
$G/H=A \oplus B$. 
Let $a\in A\backslash B$,  
$L:=\{ x\in G\mid  x+H \in B\oplus \ZZ a\}$. Set $A=\QQ a\oplus A'$, hence: 
$G/L \simeq \UU \oplus A'$. So $G/L$ embeds into $\KK$ (the restriction of this embedding 
to $U$ being identity), we equip $L$ with a regular linear order, and 
applying \cite[Lemma 5.4]{GLL 12} we get a structure 
of c-regular cyclically ordered group on $G$. The remainder of the proof is similar to the proof 
of Proposition \ref{prop53}. \end{dem}
\indent We know that an abelian group is orderable if and only if it is torsion-free, 
this is a consequence of its first order theory. 
According to \cite[Theorem 5.8]{GLL 12} being cyclically orderable 
is a consequence of the first order theory of a group. 
We are going to see that the class of infinite groups which can be equipped with a c-regular 
and discrete cyclic order and the class of infinite groups which can be equipped with a 
discrete regular linear order are not elementary.  \\[2mm]
\indent If $G$ is an infinite abelian group equipped with a structure of discrete 
c-regular cyclically ordered group,  then its torsion subgroup embeds into 
$\UU$, and $G$ contains an 
element $e$ which is not divisible by any prime $p$ and such that the 
quotient group $G/\ZZ e$ is divisible. The fact that $G$ contains an element which is 
not $p$-divisible for every prime $p$  
is not a first order property in the language of groups (the following examples 
will prove this assertion), while this is a consequence of a first order formula of 
the language of cyclically ordered group (the positive cone 
has a  smallest element). \\[2mm] 
\indent We give examples of groups which are elementarily  
equivalent to discrete c-regular cyclically ordered groups in the language of groups,  
but which cannot be equipped with a structure of discrete c-regular cyclically 
ordered group. Let $a_n$, $n\in \NN^*$, be a rationally independent family of elements of 
$\KK\backslash \UU$, and $G$ be the direct sum of the 
$\ZZ_{(p_n)}a_n$'s (where $(p_n)_{n\in \Ni^*}$ is the increasing sequence of all primes, 
and $\ZZ_{(p_n)}$ is the localization of $\ZZ$ at the ideal 
$(p_n)$). Every element of $G$ is contained in a finite sum 
$\ZZ_{(p_{n_1})}a_{n_1}+\cdots +\ZZ_{(p_{n_k})}a_{n_k}$, hence it is $p$-divisible, 
for every prime $p$ which is distinct from $p_{n_1},\dots,p_{n_k}$. 
Consequently, $G$ does not contain any element which is not divisible by any 
prime, hence it cannot be equipped with a structure of discrete c-regular
 cyclically ordered group. 
However, as we recalled in Section \ref{subsec31}, 
for every prime $p$, we have $[p]G=p$, hence $G$ is elementarily equivalent 
in the language of groups to every discrete torsion-free c-regular cyclically 
ordered group. Note that $G$ is also elementarily 
equivalent to $\ZZ$, but it cannot be equipped with a structure of discrete regular
linearly ordered group.
\begin{address}
G\'erard LELOUP\\
Laboratoire Manceau de Math\'ematiques\\
Facult\'e des Sciences\\
avenue Olivier Messiaen\\
72085 LE MANS CEDEX\\
FRANCE\\
gerard.leloup@univ-lemans.fr\\[2mm] 
Fran\c{c}ois LUCAS\\
LAREMA - UMR CNRS 6093\\
D\'epartement de Math\'ematiques\\ 
Facult\'e des Sciences\\
2 boulevard Lavoisier\\
49045 ANGERS CEDEX 01\\ 
FRANCE\\
lucasfm49@gmail.com
\end{address}
\end{document}